\documentclass[10pt]{article}
 \usepackage{amsmath, amsthm, amssymb, amsfonts}
\newtheorem{Theorem}{Theorem}

\newtheorem{Lemma}[Theorem]{Lemma}

\newtheorem{Corollary}[Theorem]{Corollary}

\newcommand{\R}{\mathbb{R}}

\newcommand{\Z}{\mathbb{Z}}
\newcommand{\I}{{\mathbf{I}}}
\newcommand{\E}{\mathbf{E}}
\newcommand{\PR}{\mathbb{P}}
\newcommand{\p}{\mathbf{p}_{\mathbf{NC}}}
\newcommand{\N}{\mathbf{N}}

\newcommand{\OH}{\hat{\Omega}}

\newcommand{\pfaff}{\mbox{Pfaff}}

\begin{document}

\title{Multi-point correlations for two dimensional coalescing random walks.}
 
\author{Jamie Lukins, Roger Tribe
and Oleg Zaboronski
\\
Mathematics Institute, University of Warwick,
Coventry, CV4 7AL, UK}

\maketitle

\begin{abstract}
This paper considers an infinite system of instantaneously coalescing rate one simple random walks on $\Z^2$,
started from the initial condition with all sites in $\Z^2$ occupied. 
Two-dimensional coalescing random walks are a `critical' model of interacting
particle systems: unlike coalescence models in dimensions three or higher, the fluctuation effects
are important for the description of large-time statistics in two dimensions, manifesting themselves
through the logarithmic corrections to the `mean field' answers. Yet the fluctuation effects
are not as strong as for the one-dimensional coalescence, in which case the fluctuation effects
modify the large time statistics at the {\it leading order}. Unfortunately,
unlike its one-dimensional counterpart,
the two-dimensional model is not exactly solvable, which explains a relative scarcity of rigorous analytic answers for
the statistics of fluctuations at large times.
This paper's contribution is to show that the correlation functions of the model
decay, for any $N \geq 2$, as
\[ \rho_N (x_1,\ldots,x_N;t) = \frac{c_0(x_1,\ldots,x_N)}{\pi^N} (\log t)^{N-{N \choose 2}} t^{-N}
\left(1 + O\left( \frac{1}{\log^{\frac12-\delta}\!t} \right) \right) 
\]
as $t \to\infty$. This generalises the results for $N=1$ due to Bramson and Griffeath
and confirms a prediction in the physics literature for $N>1$. An analogous statement holds for instantaneously annihilating random walks. 

The key tools are the known asymptotic $\rho_1(t) \sim \log t/\pi t$ due to Bramson and Griffeath, and  
the non-collision probability $\p(t)$, that no pair 
of a finite collection of $N$ two dimensional simple random walks 
meets by time $t$, whose asymptotic $\p(t) \sim c_0 (\log t)^{-{N \choose 2}}$
was found by Cox, Merle and Perkins. 
This paper re-derives the asymptotics both for $\rho_1(t)$ and $\p(t)$ by proving
that these quantities satisfy 
{\it effective 
rate equations}, that is approximate differential equations at large times.  This approach can be 
regarded as a generalisation of the Smoluchowski theory of renormalised rate equations to multi-point
statistics. 
\end{abstract}
\section{Introduction} \label{s1}
\subsection{Statement of the main result.} \label{s1.1}
We consider the following infinite coalescing systems: at time zero every site in $\Z^2$ has a particle; the particles evolve as 
independent rate one simple random walks when at disjoint sites; when any particle jumps on to another the pair 
instantaneously coalesces. 
This is a prototypical model of chemical kinetics, which has been studied intensively by both the physics and
mathematics communities, see e.g. \cite{tauber} and \cite{CG} for reviews. Consequently, a great deal is already known
about the statistical behaviour of coalescence in all dimensions. The only steady state of the model
is the absorbing state with zero particle density, and the {\it approach} to this trivial state is a first interesting question about coalescence.
In particular, the particle density $\rho_1(x;t) = \PR[ \mbox{$x$ is occupied at time $t$}]$ is independent of
$x$ for a fully occupied initial condition, and satisfies in dimension two
\begin{equation} \label{sawyer}
\rho_1(t) =  \frac{\log t}{\pi t}  \left(1 + O\left( \frac{1}{\log^{\frac12}\!t} \right) \right) \qquad \mbox{as $t \to\infty$.}
\end{equation}
The leading asymptotic in (\ref{sawyer}) is Theorem $1^{\prime}$ of Bramson and Griffeath \cite{Bramson+G} who 
derive it using asymptotics of Sawyer \cite{Sawyer} for the voter model. 
Note the structure of the leading term: within the standard assumptions of chemical kinetics which disregard correlations between reactants (the `well-mixed' limit, see
e.g. \cite{chembook}),
the evolution of $\rho_1(t)$ is described by the `mean field rate equation' $$\dot{\rho}_1(t)=-\lambda \rho_1^2(t),$$
where $\lambda>0$ is an effective reaction rate. The mean field prediction is therefore
that $\rho_1(t)$ decays as $t^{-1}$ in the large time limit. The true answer (\ref{sawyer}) agrees
with this prediction up to a logarithmic correction only, which tells us that correlation effects are important in two
dimensions. To quantify these effects in more detail,
we study the large-time behaviour of multi-point correlation functions 
\[
\rho_N (x_1,\ldots,x_N;t) := \PR[ \mbox{$x_1,\ldots,x_N$ are occupied at time $t$ }\!\!].
\]
The following statement is the main result of this paper:
\begin{Theorem} \label{T1}
For $N \geq 2$ and disjoint $x_1,\ldots,x_N$, there exists $c_0(x_1,\ldots,x_N)$ so that for any $\delta \in (0,\frac12)$,
as $t \to \infty$
\begin{eqnarray} \label{MR2}
\rho_N (x_1,\ldots,x_N;t) & = &  c_0(x_1,\ldots,x_N)  \, \rho_1^N(t) (\log t)^{-{N \choose 2}} 
\left(1 + O\left( \frac{1}{\log^{\frac12- \delta}\!t} \right) \right) \nonumber \\
& = & \frac{c_0(x_1,\ldots,x_N)}{\pi^N} (\log t)^{N-{N \choose 2}} t^{-N}
\left(1 + O\left( \frac{1}{\log^{\frac12-\delta}\!t} \right) \right).
\end{eqnarray}
\end{Theorem}
The well known thinning relation between coalescing and annihilating systems (reviewed for example in  \cite{TZ1} Section 2.1) 
allows one to apply Theorem \ref{T1} to the study of
an instantaneously annihilating system. Under this thinning, the initial condition for a coalescing particles system 
of all sites occupied corresponds to an annihilating system with a Bernouilli ($\frac12$) initial condition. 
Theorem \ref{T1} then immediately implies the following statement:
\begin{Corollary} \label{T1a}
Consider the following infinite annihilating system: at time zero every site in $\Z^2$ is independently 
occupied with probability $\frac{1}{2}$; the particles evolve as independent rate one simple random 
walks when at disjoint sites; when any particle jumps on to another the pair 
instantaneously annihilates (disappears from the system). Let $\rho_N^a (x_1,\ldots,x_N;t)$ be the probability that sites
$x_1, x_2, \ldots, x_N$ are occupied at time $t$. 
For $N \geq 2$, disjoint $x_1,\ldots,x_N$ and for any $\delta \in (0,\frac12)$,
as $t \to \infty$
\begin{eqnarray} \label{MR2+}
\rho_N^a (x_1,\ldots,x_N;t) =
\frac{c_0(x_1,\ldots,x_N)}{(2\pi)^N} (\log t)^{N-{N \choose 2}} t^{-N}
\left(1 + O\left( \frac{1}{\log^{\frac12-\delta}\!t} \right) \right).
\end{eqnarray}
\end{Corollary}
\noindent
{\bf Remark.} Alternatively, to treat the fully occupied initial condition for the annihilating system, one can repeat 
the proofs in this paper. The effective rate equation for $\rho_1$ has an extra power of $2$, namely (\ref{B1}) is replaced by 
$\frac{d}{dt} \rho_1(t) = - 2 \E [ \xi_t(0) \xi_t(e_1)]$. This extra factor of $2$ is then the only change and traces through
to become a factor of  $2^{-N}$  by the time one reaches the asymptotics for $\rho_N$. 

The main idea of the proof of Theorem \ref{T1}, see the start of Section \ref{s1.2}, is to reduce the question about an infinite coalescence 
system to the one concerning finitely many non-interacting particles, as follows:  for $N\geq 2$ 
simple random walks on $\Z^2$ started from fixed 
distinct positions $x_1, x_2, \ldots, x_N \in \Z^2$ estimate the non-collision probability
\[
\p(s) = \PR[ \mbox{no pair of walkers meets by time $s$}]. 
\]
The key ingredient for Theorem \ref{T1} is the following asymptotic for $\p$, which shows exactly the
source of the constant $c_0(x_1,\ldots, x_N)$ in Theorem \ref{T1}, and is due to Cox, Merle and Perkins (see \cite{CMP} Prop. 1.3.):
there exists a positive constant $c_0(x_1,\ldots, x_N)$ so that
\begin{equation} \label{MR1}
\p(t) =  c_0(x_1,\ldots, x_N) (\log t)^{-{N \choose 2}} \left(1 + O\left( \frac{\log \log t}{\log t} \right) \right) \qquad \mbox{as $t \to\infty$.}
\end{equation}
This is a foundational result about two-dimensional random walks which has 
applications far beyond the interacting particle systems
considered in the present paper and the Lotka-Volterra model studied in \cite{CMP}. 
We discovered that both the asymptotic for $\rho_1$ in (\ref{sawyer}) and the non collision probability in
(\ref{MR1}) above, follow using the formalism of effective rate equations. We present these new derivations,
partly as a way to advertise (\ref{MR1}), which is currently rather hidden inside \cite{CMP}. 

\noindent
{\bf Paper organisation.}  In section \ref{New} we discuss the main result, the non-collision probabilities, and their
companions in dimensions $d \geq 3$ and $d=1$, explaining why the case of $d=2$ is comparatively harder to analyse.
Section \ref{s1.2} outlines the main idea behind proof of Theorem \ref{T1}, and of 
the asymptotics (\ref{sawyer}) and (\ref{MR1}). Section \ref{s2} contains detailed error estimates 
that complete the proofs. 

\subsection{Discussion and related results.} \label{New}

Theorem \ref{T1} shows that the multipoint correlation $\rho_N$ decreases with time faster than $\rho_1^N$, which
reflects the negative correlation for occupation numbers of nearby sites. This negative correlation is a well
known property of coalescing random walks in any number of dimensions, see e.g. \cite{KvdB1}. 
In addition, equation (\ref{MR2}) quantifies the nature of negative correlations
in two dimensions: the fact that $\rho_N/\rho_1^N$ decays as 
$(\log t)^{-{N \choose 2}} $, suggests that negative correlations in two dimensions have a pairwise structure, meaning
that the effect of particles $A$ and $B$ on particle $C$ can be correctly reproduced at the leading order by 
disregarding correlations between $A$ and $B$.

It is worth noting, that the statement of Theorem \ref{T1} is hard to reproduce employing traditional phenomenological
theories used in physics to capture correlation effects in interacting particle systems. The most famous example
of such a phenomenology is Smoluchowski's approximation reviewed in \cite{KBNR}, which states that the mean field approximation
$\rho_N=\rho_1^N$ in the right hand side of the Hopf equation for $\rho_N$ 
should be modified as follows: $\rho_N=N\lambda_{eff} (t) \rho_1^N$, where $\lambda_{eff}(t)$
is the effective reaction rate equal to the flux of Brownian particles seeded with density $\rho_1$
at infinity through a sphere of a fixed radius $r_0$ around the target particle.
In Smoluchowski theory, $r_0$ is interpreted as the effective interaction radius.
 In two dimensions, $\lambda_{eff}(t)=\frac{2\pi \rho_1}{\log(\sqrt{t}/r_0)}$. Solving the resulting equation
 for $\rho_1$, one finds that $\rho_1\sim \frac{\log(t)}{t}$, which gives the correct law of decay (\ref{sawyer})
 albeit with a wrong constant factor in the leading term. However, for $N>1$, the solution of the
 rate equations modified according to Smoluchowski's theory gives $\rho_N=\rho_1^N$, which contradicts
 (\ref{MR2}).

It was argued in \cite{RZ} that the fluctuation effects for models in the same universality class as coalescing random 
walks in one and two
dimensions are encoded not just in the effective reaction rate, but crucially in the renormalisation
of composite field operators which describe occupation probabilities in the limit of large times using the language
of effective quantum field theory. This paper uses a non-rigorous method of dynamical renormalisation
group which accounts both for the emergence of Smoluchowski's reaction rate and the anomalous
dimensions of the composite field operators to predict the correct decay law (\ref{MR2}).  
One of the aims of the present work is thus to provide a rigorous proof of (\ref{MR2}).

The study of multi-site occupation probabilities is uniquely difficult in two-dimensions: in dimensions three and
higher it is widely expected that fluctuation effects do not invalidate the (properly understood) 
mean field assumptions, consequently $\rho_N$ decays with time as $\rho_1^N$, see Section \ref{s1.3}
for more details. In one dimension the fluctuation effects are important, but the model is exactly solvable
leading to exact expressions for multi-site probabilities in terms of Pfaffians, see \cite{TZ2}. The derivation of
the counterpart of (\ref{MR2}) in one dimension thus reduces to an exercise in asymptotic analysis.
For the flavour of the one-dimensional computations see below.

 To account for fluctuation effects in two dimensions we adapt the method of van der Berg and Kesten \cite{KvdB1}
which was originally developed for dealing with interacting particle systems in dimension higher than $2$.
Therefore, it is not unreasonable to predict that the method we use below should work for $d \geq 3$ that
$\rho_N (x_1,\ldots,x_N;t) \sim c_0(x_1,\ldots,x_N)  \, \rho_1^N(t)$. What is perhaps more surprising,
 similar methods can also applied to one-dimensional coalescence leading to upper and lower 
 bounds with correct time dependence, albeit with constants not close to the exact value, see \cite{TZ1} for details.

To conclude the discussion of our results, let us stress that the question about the spatial dependence
of the occupation probabilities is still very much open. We still do not know the full behaviour of 
$c_0(x_1,\ldots,x_N)$ as a function of $(x_i)$'s , and it
is perhaps surprising that the large time asymptotics can be found without any such knowledge. 
A reasonable conjecture however, consistent with the results for well separated initial conditions and the upper bound found in 
\cite{CMP}, would be that
\begin{equation} \label{conj}
c_0(x_1,\ldots,x_N) \approx \prod_{i<j} \log(|x_i-x_j|^2) \quad \mbox{for $O(1) << |x_i-x_j| << O(t^{1/2}).$}
\end{equation}
Notice that such a product form ansatz was already used in \cite{Gaudilliere} to construct a sub-harmonic function
on a subset of $\R^{2N}$ corresponding to non-collision for the analogous problem of collision of Brownian discs. 
Notice also that the structure of the guess (\ref{conj}) is similar to the $d=1$ result (\ref{oned}), see below.
The crucial modification contained in (\ref{conj}) in comparison with \cite{Gaudilliere} 
is that the logarithmic dependence of $\p$
on the separation between points is $asymptotic$ and does not 
extend all the way to the boundary of the non-collision region.

The non-collision probability $\p$ has also been discussed for random walks in all dimensions $\Z^d$, and for the closely related problems of
Brownian motions in $d=1$ or Brownian discs in $d>1$,  but the argument strongly depends on $d$.
In dimension $d \geq 3$, the problem is easy since one has $\lim_{t \to \infty} \p(t) >0$ due to the transience of the underlying walks. 
Indeed, letting $\Omega_{ij}$ be the event that particles $i$ and $j$ ever meet, 
for sufficiently widely separated starting points a simple union bound shows that $\PR[\cup_{i<j} \Omega_{ij}] <1$.
Starting from this remark, a lot can be said about $\p$ using standard probabilistic tools. In particular, the statement
that $\rho_N (x_1,\ldots,x_N;t)  \sim \rho_1^N $ can be established.

In dimension $d=1$ the correlation effects are significant, but the problem is 'exactly solvable' and can be fully analysed using the celebrated Karlin-McGregor formula. The asymptotics for
the Brownian analogue of $\p(t)$ is discussed in Grabiner  \cite{grabiner}.  
To give the flavour of the `integrable probability' , we will review the argument here, but streamline it slightly by employing the De Bruijn integration formula 
\cite{debruijn}. T
To avoid unnecessary technicalities, let us assume that $N$ is even.
Recall that the Karlin-McGregor formula states that the transition density for Brownian particles in $\R$
started at $x_1<x_2<\ldots <x_N$ to end up at $y_1<y_2<\dots<y_N$ at time $t$ without
meeting each other is given by $\det_{1\leq i,j \leq N} [g_t(x_i, y_j)]$, where $g_t$ is
the one-dimensional Gaussian transition density. Therefore,
\begin{eqnarray}\label{onedexact}
\p(t)&=&\int_{-\infty<y_1< y_2<\ldots < y_N<\infty} dy_1 dy_2 \ldots dy_N
\det_{1\leq i,j \leq N} [g_t(x_i, y_j)]\nonumber
\\&=&\pfaff_{1\leq i<j\leq N}\left[\Phi\left(\frac{x_i-x_j}{\sqrt{4t}}\right)\right],
\end{eqnarray}
where $\Phi(x)=2\int_0^x \frac{dt}{\sqrt{\pi}}e^{-t^2}$ and $\pfaff(M)$ is the Pfaffian
of an anti-symmetric matrix $M$. The second equality follows from
the De Bruijn formula. 
This is an 
exact answer for the non-colliding probability, which is  well suited for the asymptotic analysis.
The individual terms in the Laplace expansion of the Pfaffian decay as $t^{-N/2}$, so a lot of cancellation should occur before one sees the correct decay. 
However it is immediately seen from (\ref{onedexact}) that $\p(\cdot,t)\mid_{x_i=x_{j}}=0$,
suggesting that the Vandermonde determinant $\Delta(x_1,\ldots,x_N) := \prod_{i<j}(x_i-x_j)$ would factor out of
a power series expansion around the origin.  A straightforward analysis (exploiting a suitable conjugation of the matrix 
$M$ - see \cite{TZ2} Lemma $5$) 
shows that
\begin{eqnarray}\label{oned}
\p(x_1,\ldots, x_N,t)&=&c_N \Delta\left(\frac{x_1}{\sqrt{t}}, \frac{x_2}{\sqrt{t}},\ldots, \frac{x_N}{\sqrt{t}}\right)\left(1+O\left(\frac{x_{ij}}{\sqrt{t}}, 1\leq i<j\leq N\right)\right),
\end{eqnarray}
where $c_N$ is a constant
the precise value of which is known \cite{grabiner}.
So the non-collision probability decays as $t^{-{N \choose 2}}$ at large times.
Notice that the proportionality coefficient $\Delta(x_1,\ldots,x_N)$ 
is a harmonic function on the  Weyl chamber $W_N=\{-\infty<x_1< x_2<\ldots < x_N<\infty\}$
equal to zero at the boundary $x_i=x_{i+1},~1\leq i<N$, corresponding to collisions. 
It is interesting to point out that a (false) assumption 
that pairwise collisions occur independently of each other leads to the correct decay rate in (\ref{oned}) but not the
correct constant $c_N$.

For $d=2$, there are no counterparts of Karlin-McGregor and DeBruijn formulae. Simple union bounds
which work for $d>2$ also fail to capture the correct asymptotic of $\p$ due to the relevance of correlations.  
This probably explains why the two-dimensional question was solved only recently in \cite{CMP}, where it appears
as key tool  for the study of the class of models which can be regarded as perturbations of a
two dimensional voter model. For the sake of completeness, let us also mention two related results: for random
walks with initial separations growing as a power law of time, the asymptotic independence holds and $\p(t)$
becomes soemwhat easier to estimate,  see the survey by Cox and Griffeath \cite{CG}.
The correct decay of $\p$ has also been conjectured in \cite{Gaudilliere},
but the author only managed to prove a weaker lower bound $\p(t)>t^{-cN^4}$ for some positive $c$. 


\subsection{Sketch of the main arguments.} \label{s1.2}

The decay rate of the particle intensity $\rho_1(t)$ was re-examined by Kesten and van den Berg \cite{KvdB1}, \cite{KvdB2} in dimensions $d \geq 3$. They 
showed that the derivative $d \rho_1(t)/dt$ can be approximated at large times so that $\rho_1(t)$ approximately solves an
autonomous ODE, from which the asymptotics can be derived. 
We show below that this effective rate equation approach also works also in dimension $d=2$, to yield the 
result (\ref{sawyer}). 
We think this alternative argument worth recording, not just for its intuitive nature, but because (as shown by Kesten and 
van den Berg in dimensions $d \geq 3$) it has the potential to apply to related systems where the original duality methods do not apply. 
However we consider only the simplest case (instantaneous simple random walks)
we can apply the negative correlation properties discussed above, 
which reduce the trickier variance estimate needed in \cite{KvdB1} to a triviality using negative correlation properties.
Moreover the leading rate and leading constant in the asymptotic (\ref{sawyer}) for $d=2$ has a degree of universality 
in that, due to the recurrence, they should not depend on the exact nature of the coalescence mechanism, 
and on only the variance of the underlying walk, not the exact choice of step distribution or underlying lattice 
which only affect sub-leading terms in the asymptotic. Some of these questions are explored in the \cite{thesis}, 
and furthermore we leave some details of the current sketch to this upcoming thesis. 

To our surprise the effective rate equation also provides a natural approach to the asymptotic (\ref{MR1}) for the non-collision
probability.  We include a proof using an effective rate equation 
as an alternative to the clever repeated conditioning argument used in \cite{CMP}.

We give the main heuristics behind both effective rate equations below, leaving the estimates on the 
various error terms $\mathcal{E}_1,\mathcal{E}_2, \ldots$ to section \ref{s2}. We note that the method allows one easily to track error estimates,
which we have included in the statements (\ref{sawyer}) and (\ref{MR1}) (without believing in some cases that they are best possible). 

\noindent
\textbf{Effective rate equation for $\p$.}
We list the positions of $N$ independent simple random walks as $((X^{(i)}_t): i = 1,\ldots,N)$.
We suppose $X^{(i)}_0=x_i$ for $i \leq N$ under the probability measure $\PR_{x}$, where $x=(x_1,\ldots,x_N)$. 
We write $\tau_{ij}$ for the collision time between $X^{(i)}$ and $X^{(j)}$; we write
$\tau_c = \min_{i \neq j} \tau_{ij}$ for the first collision time between any pair. 

Using the generator for the walk $X=(X^{(1)}, \ldots,X^{(N)})$ killed upon the first collision we find
\begin{equation} \label{gencalc}
\frac{d\p}{dt}(t)  = - \frac12 \sum_{i<j} \PR_x[\tau_{c}>t, |X^{(i)}_t-X^{(j)}_t|=1].
\end{equation}
We introduce an intermediate time $s \in [0,t]$, which we choose carefully later
(in fact we end up choosing $s = t  \log^{-\alpha}\!t$). 
A-posteriori we know that 
$\p(s) \approx \p(t)$, which motivates the next approximation 
\[
\frac{d\p}{dt}(t)  \approx - \frac12 \sum_{i<j} \PR_x[\tau_{c}>s, \tau_{ij}>t, |X^{(i)}_t-X^{(j)}_t|=1]
\]
(we use the symbol $\approx$ informally to mean approximately equal). 
The Markov property at time $s$  yields
\[
\frac{d\p}{dt}(t)  \approx - \frac12 \sum_{i<j} \E_x [\I(\tau_c>s) F_{2(t-s)}(X^{(i)}_s-X^{(j)}_s)] 
\]
where 
\[
F_t(y) = \PR_y[|X_t|=1,\tau_0>t]
\]
 for a single rate one walker on $\Z^2$ and its hitting time $\tau_0$ of the origin.
An analysis of $F$ will show that $F_{2(t-s)}(y) \approx \frac{2}{t \log t} $ for most of the values
$y = X^{(i)}_s-X^{(j)}_s$ (which relies on the choice $s = t \log^{-\alpha}\!t$). 
This will decouple the event $\{\tau_c >s\}$ from the variable $F_{2(t-s)}(X^{(i)}_s-X^{(j)}_s)$ and allows a further
approximation
\[
\frac{d\p}{dt}(t)  \approx - {N \choose 2}  \PR_x [\tau_c>s]  \frac{1}{t \log t} = - {N \choose 2}  \frac{1}{t \log t} 
\p \left(t \log^{-\alpha}\!t\right). 
\]
The details in the section \ref{s2.1} quantify the approximations above and establish the following 
approximate rate equation
\begin{equation} \label{ERE1}
\frac{d\p}{dt}(t)  = - {N \choose 2} \frac{1}{t \log t} \p \left(t \log^{-\alpha}\!t\right) (1+\mathcal{E}_1(t)) + \mathcal{E}_2(t)
\end{equation}
where the errors $\mathcal{E}_i$ are bounded as follows
\begin{equation} \label{errorbounds}
\mathcal{E}_1(t)  = O\left( \frac{ \log \log t}{\log t}\right)  \qquad \mathcal{E}_2(t) = O \left(\frac{1}{t \log^{\mu}\!t}\right) \quad \mbox{for any $\mu>0$}.
\end{equation}
It is now a calculus exercise  to show that (\ref{ERE1})
is sufficiently close to the effective rate equation
\[
\frac{d\p}{dt}(t)  = - {N \choose 2} \frac{1}{t \log t} \, \p(t)
\]
that the desired asymptotic (\ref{MR1}) and associated error bound for $\p$ can be deduced.

\noindent
\textbf{Effective rate equation for $\rho_1$.}

We mimic closely the heuristic argument from Kesten and van den Berg \cite{KvdB1}.
The initial condition implies that $\rho_1(x;t)$ is independent of $x$. We write $x \sim y$ to mean that two sites $x,y \in \Z^2$ are neighbours. 
An exact generator calculation shows (writing $e_1 =(1.0)$)
\begin{eqnarray*}
\frac{d}{dt} \rho_1(0;t) &=&  - \frac14 \sum_{x \sim 0} \PR[ \mbox{sites $0$ and $x$ are occupied at time $t$}] \\
& = & - \PR[ \mbox{sites $0$ and $e_1$ are occupied at time $t$}].
\end{eqnarray*}
To find the effective rate equation on traces back the paths of particles found at $0$ and $e_1$ at time $t-s$
to an earlier time $s$, noting that over $[t-s,t]$ they must not coalesce. There may be more than one particle
at time $t-s$ that end up at the origin at time $t$, but choosing $s$ suitably there are likely to be only one, and we reach the approximation
\[
\frac{d}{dt} \rho_1(0;t) \approx - \sum_{x \neq y} \PR[\mbox{sites $x$ and $y$ are occupied at time $t-s$}] \psi_s(x,y)
\]
where
\[
\psi_s(x,y) = \PR_{(x,y)}[ X^{(1)}_s = 0, X^{(2)}_s = e_1, \tau >s] = \PR_{(0,e_1)}[ X^{(1)}_s = x, X^{(2)}_s = y, \tau >s] 
\]
The main difference in dimension $d=2$ is to justify the approximations
\[
\psi_s(x,y) \approx p_s(x) p_s(y) \PR_{(0,e_1)}[\tau >s] \approx 
p_s(x) p_s(y) \frac{\pi}{\log s}
\]
for typically located $x,y$ (that is not too close) and 
\[
\PR[\mbox{sites $x$ and $y$ are occupied at time $t-s$}] = \rho_2(x,y;t-s) \approx \rho_1(x;t-s) \rho_1(y;t-s)
\]
again for suitably well spaced $x,y$. These approximations lead to 
\[
\frac{d}{dt} \rho_1(0;t) \approx - \frac{\pi}{\log s} \rho_1(t-s)^2.
\]
Choosing the intermediate time $s = t \log^{-\frac12} \!t$ the error bounds  in section \ref{s2.2} will show that 
\[
\frac{d}{dt} \rho_1(t) =   - \frac{\pi}{\log t} \rho_1^2(t) + O\left(\frac{\log^{\frac12}\!t}{t^2}\right)
\]
and that this implies the desired asymptotic (\ref{sawyer}) and associated error bound. 

\noindent
\textbf{Application to the decay rate of $\rho_N$.}

When $x_1, \ldots,x_N$ are occupied at time $t$ we may trace paths of these particles 
backwards to an earlier time $t-s$, over which interval they must not coalesce. Although for a coalescing system
these backwards paths are not unique, by choosing the earlier time carefully we find (by an inclusion exclusion argument)
the paths are with high probability unique. This leads to the approximation
\begin{eqnarray*}
\rho_N(x_1,\ldots,x_N;t) & \approx & \sum_{y_i,\ldots,y_N \in \Z^2} \rho_N(y_1,\ldots,y_N;t-s) 
P_{(y_1,\ldots,y_N)}[ X^{(i)}_s =x_i, \tau^c > s] \\
 & = & \sum_{y_i,\ldots,y_N \in \Z^2} \rho_N(y_1,\ldots,y_N;t-s) 
P_{(x_1,\ldots,x_N)}[ X^{(i)}_s =y_i, \tau^c > s]
\end{eqnarray*}
In this sum we will justify replacing $ \rho_N(y_1,\ldots,y_N;t-s)$ by $\rho_1^N(t-s)$, using the fact the sum is dominated by
$(y_1,\ldots,y_N)$ with coordinates spaced at a distance $O(s^{1/2})$. This leads immediately to 
\[
\rho_N(x_1,\ldots,x_N;t) \approx \rho_1^N(t-s) \p(s)
\]
and with the suitable choice of the intermediate time $t-s$ the result follows from the asymptotics for
$\rho_1$ and $\p$.

\section{Details.} \label{s2}
We justify the effective rate equations for $\p$ in section \ref{s2.1}, and for $\rho_1$ in section \ref{s2.2}.
We give the proof of Theorem \ref{T1} in section \ref{s2.3}.

\noindent
\textbf{Notation.}
We write $p_t(y)$ for the transition probabilities of a single simple random walk on $\Z^2$. When $y=(y_1,\ldots,y_K)$
we also write $p_t(y)$ for $\prod_{i \leq K}p_t(y_i)$ the transition probabilities for a system of independent walks.

\noindent
We write $\mathcal{E}_1,\mathcal{E}_2,\ldots$ for various error terms that arise in the paper; the subscripts are restarted
in each subsection. 

\noindent
Throughout $c_0,c_1, \ldots$ are constants to which we may refer, whereas $C$ is a running constant that may 
change value from line to line. We indicate the dependencies of constants upon parameters, \textit{except} that $N \geq 1$ and 
the positions $x=(x_1,\ldots,x_N)$ are fixed throughout and we often suppress notation indicating dependency on these quantities. 

\subsection{Non-collision for finite particle system} \label{s2.1}
We start by completing the calculus exercise for the effective rate equation 
assuming the error bounds (\ref{errorbounds}). Choose $\mu = {N \choose 2} +2$ 
and choose $t_0 \geq 2$ (ensuring $\log t_0>0$) so that $1+|\mathcal{E}_1(t)| \geq \frac12$ for all $t \geq t_0$. 
Then using $\p(t \log^{-\alpha}\!t) \geq \p(t)$ we obtain the upper bound
\[
 \frac{d\p}{dt}(t)  \leq - {N \choose 2}  \frac{\p(t)}{t \log t} (1+\mathcal{E}_1(t)) + \mathcal{E}_2(t) \quad \mbox{for $t \geq t_0$.}
\]
Using the integrating factor $\exp(I(t))$ where
\[
I(t) =  {N \choose 2} \int^t_{t_0} \frac{1+\mathcal{E}_1(s)}{s \log s} ds 
\]
we obtain 
\[
\p(t) \leq e^{-I(t)} \left( \p(t_0) e^{I(t_0)} +  \int^t_{t_0} e^{I(s)} \mathcal{E}_2(s) ds \right).
\]
Using the bounds (\ref{errorbounds}) on $\mathcal{E}_i$ we see that 
$ I(t) = {N \choose 2} \log \log t + O(1)$ and that $\int^t_{t_0} e^{I(s)} |\mathcal{E}_2(s)| ds $ is bounded as $t \to \infty$. 
This implies the existence of $C$ so that
\[
\p(t) \leq C (\log t)^{-{N \choose 2}} \qquad \mbox{for all $t \geq t_0$.}
\] 
We now use this upper bound to help estimate the error in replacing $\p(s)$ by $\p(t)$ when
$s=t \log^{-\alpha}\!t$.
Using this upper bound and (\ref{ERE1}) we see that for all $t$ sufficiently large
\[
\left| - \frac{d\p}{dt}(t) \right|  \leq  C \frac{1}{t (\log t)^{{N \choose 2}+1}} 
\]
so that
\begin{eqnarray*}
|\p(s) - \p(t)| & = &   \int^t_s - \frac{d\p}{dr} dr \\
& \leq & C \int^t_s \frac{1}{r (\log r)^{{N \choose 2}+1}} dr \\
& = & C \left( (\log s)^{-{N \choose 2}}-  (\log t)^{-{N \choose 2}}\right) \\
& \leq & C \frac{\log \log t}{(\log t)^{{N \choose 2}+1}}.
\end{eqnarray*}
This allows us to replace $\p(s)$ by $\p(t)$ in (\ref{ERE1}) to reach
\begin{equation} \label{ERE22}
\frac{d\p}{dt}(t)  = - {N \choose 2} \p(t) \frac{1}{t \log t}(1+\mathcal{E}_1(t)) + \mathcal{E}_3(t)
\end{equation}
where 
\[
 \mathcal{E}_1(t) = O \left(\frac{\log \log t}{\log t} \right)
 \qquad
 \mathcal{E}_3(t) = O \left( \frac{\log \log t}{t (\log t)^{{N \choose 2}+2}}\right).
 \]
We now repeat the integrating factor argument, taking a bit more care. First 
\[
I(t)= {N \choose 2} \log \log t + \mathcal{E}_4(t) 
\quad \mbox{where} \quad
\mathcal{E}_4(t) = {N \choose 2}  \int^t_{t_0} \frac{\mathcal{E}_1(s)}{s \log s} ds. 
\]
The limit $\mathcal{E}_4(\infty) = \lim_{t \to \infty} \mathcal{E}_4(t)$ exists and
$|\mathcal{E}_4(\infty) - \mathcal{E}_4(t)| = O(\log \log t/\log t)$. Then
\[
e^{-I(t)} = e^{-\mathcal{E}_4(\infty)} (\log t)^{-{N \choose 2}} \left(1 + O\left( \frac{\log \log t}{\log t} \right) \right).
\]
The integrating factor gives
\[
\p(t) =  e^{-I(t)} \left(p(t_0) e^{I(t_0)} + \mathcal{E}_5(t) \right)
\quad 
\mbox{where} \quad 
\mathcal{E}_5(t) = \int^t_{t_0} e^{I(s)} \mathcal{E}_3(s) ds.
\]
Again $\mathcal{E}_5(\infty) = \lim_{t \to \infty} \mathcal{E}_5(t)$ exists and
$|\mathcal{E}_5(\infty) - \mathcal{E}_5(t)| = O(\log \log t/\log t)$. 
Now this can be rearranged to see that 
\[
\p(t) = c_0 (\log t)^{-{N \choose 2}} \left(1 + O\left( \frac{\log \log t}{\log t} \right) \right) 
\]
where $c_0 = \left(\p(t_0) e^{I(t_0)} + \mathcal{E}_5(\infty)\right) e^{-\mathcal{E}_4(\infty)}$, completing the proof
of (\ref{MR1}).

To establish the error bounds on $\mathcal{E}_1(t), \mathcal{E}_2(t)$ in the effective rate equation we need 
estimates on the killed transition 
 density for a single random walker which we state here.  
The proof is delayed until the end of this subsection. 
\begin{Lemma} \label{Gbounds}
Let $\tau_0$ be the hitting time of the origin for $(X_t)$ a rate one simple random walk on $\Z^2$ started at $y$ under the probability measure $\PR_y$. 
Recall we write $p_t(y)$ for the transition probability.
Let $q_t(y,z) = \PR_y[X_t=z, \tau_0 >t]$ be the transition probability for the process killed at $\{0\}$. 
\begin{itemize}
\item[(a)] For $e_1=(1,0)$ and $e_2=(0,1)$ we have 
\[
\PR_{\pm e_i}[\tau_0 >t] = \frac{\pi}{\log t} + O(\log^{-2}\!t)   \quad \mbox{as $t \to \infty$.}
\]
Also for any $\beta>0$ there exists $c_1 = c_1(\beta)$ so that 
\[
\PR_{y}[\tau_0 \leq t]  \leq \frac{c_1 \log \log t}{\log t} \qquad 
\mbox{for all $t \geq 4$ and $\frac{t}{\log^{\beta}\!t} \leq |y|^2 $.}
\]
\item[(b)]  There exists $c_2,c_3$ so that
\[
q_t(\pm e_i,z) \leq \frac{c_2 p_t(z)}{\log t} + \frac{c_2}{t \log^2 \! t} 
\leq \frac{c_3}{t \log t} 
\qquad \mbox{for all $t \geq 4$ and $z$.}
\]
\item[(c)] For any $\beta_1 > \beta_2 \geq 1$ there exists $c_4= c_4(\beta_1,\beta_2)$ so that
\[
\left| q_t(\pm e_i,z) - \frac{1}{t \log t} \right|  \leq \frac{c_4 \log \log t}{t \log^2 \! t} 
\qquad \mbox{for all $t \geq 4$ and $\frac{t}{\log^{\beta_1}\!t} \leq |z|^2 \leq \frac{t}{\log^{\beta_2}\!t}$.}
\]
\end{itemize}
\end{Lemma}

We now establish the effective rate equation (\ref{ERE1}).
For convenience we use the notation
$ X^{(i,j)}_s = X^{(i)}_s - X^{(j)}_s$ (and $x_{ij} = x_i - x_j$).
We repeat the main argument - there are three approximation steps:
\begin{eqnarray}
\frac{d}{dt} \p(t) & = &  - \frac12 \sum_{i<j} \PR_x[\tau_{c}>t, |X^{(i,j)}_t|=1] \nonumber \\
& \approx & - \frac12 \sum_{i<j} \PR_x[\tau_{c}>s, \tau_{ij}>t, |X^{(i,j)}_t|=1]  \label{app1} \\
& = &  - \frac12 \sum_{i<j} \E_x [\I(\tau_c>s) F_{2(t-s)}(X^{(i,j)}_s)]\nonumber \\
& \approx &  - \frac12 {N \choose 2}  \PR_x [\tau_c>s]  \frac{4}{2(t-s) \log 2(t-s)} \label{app2} \\
& \approx & - {N \choose 2} \frac{1}{t \log t} \p(t \log^{-\alpha}\!t) .\label{app3}
\end{eqnarray}
We now derive the error bounds (\ref{errorbounds}) for these approximations, where additive errors are collected in $\mathcal{E}_2$
and multiplicative errors are collected in $(1+\mathcal{E}_1)$, as written in (\ref{ERE1}). We set $s = t \log^{-\alpha}\!t$ (fixing
$\alpha>0$ soon) and we consider only $t$ where $0 < s < \frac12 t$.

\noindent
\textbf{Step 1.} For the error in the first approximation (\ref{app1}), for any $i<j$ we may bound 
\begin{eqnarray}
0 & \leq & \PR_x[\tau_{c}>s, \tau_{ij}>t, |X^{(i,j)}_t|=1]  
 - \PR_x[\tau_{c}>t, |X^{(i,j)}_t|=1]  \nonumber \\
& \leq &
\sum_{k<l}  \PR_x[ \tau_{c}>s, \tau_{ij}>t, |X^{(i,j)}_t|=1, \tau_{kl} \leq t].  \label{E1}
\end{eqnarray}
We consider first the case where $\{k,l\} \cap \{i,j\} = \emptyset$.
Then,  by the Markov property at time $s$,
\[
\PR_x[\tau_c>s, \tau_{ij}>t, |X^{(i,j)}_t|=1, \tau_{kl} \leq t] \\
= \E_x[ \I(\tau_c>s) F_{2(t-s)}(X^{(i,j)}_s) G_{2(t-s)}(X^{(k,l)}_s)]
\]
where $G_t(z) = \PR_z[ \tau_0 \leq t]$,
$q_t(y,z) = \PR_y[X_t=z, \tau_0 >t]$ is the transition probability killed at $\{0\}$, and 
\begin{equation} \label{Fdefn}
F_t(y) = \PR_y[|X_t|=1,\tau_0>t] = q_t(y,e_1) + q_t(y,-e_1) + q_t(y,e_2) + q_t(y,-e_2). 
\end{equation}
Applying the bounds from Lemma \ref{Gbounds} (b) we find 
\begin{eqnarray*}
&& \hspace{-.3in} \PR_x[\tau_{c}>s, \tau_{ij}>t, |X^{(i,j)}_t|=1, \tau_{kl} \leq t] \\
  & \leq & \frac{4c_3}{2(t-s) \log 2(t-s)} \E_x[\I(\tau_c>s) G_{2(t-s)}(X^{(k,l)}_s)] \\
& \leq & \frac{C}{t \log t} \E_x[\I(\tau_c>s) G_{2(t-s)}(X^{(k,l)}_s)] 
\end{eqnarray*}
We break this expectation into two parts: 
\[
\E_x\left[ \I(\tau_c>s) G_{2(t-s)}(X^{(k,l)}_s) \I( s^{-1/2} X^{(k,l)}_s \in S_1 \cup S_2) \right]
\]
where
\[
S_1 = \{|z| \leq \log^{-\gamma}\!s \}, \quad
S_2 = \{|z| > \log^{-\gamma}\!s \}.
\]
For the part $S_1$, by a central limit theorem estimate (for example see (\ref{lclt}))
\[
\PR_x[s^{-1/2} X^{(k,l)}_s \in S_1] \leq C \log^{-2\gamma} \!s
\]
 so that by choosing $2 \gamma +1 = \mu$ this part contributes to $\mathcal{E}_2$. 
For the part corresponding to $S_2$, we choose $\beta > 2\gamma + \alpha$ and the estimate on
$G$ in Lemma \ref{Gbounds} (a) applies for large $t$, giving
\begin{eqnarray*}
\E_x[ \I(\tau_c>s) G_{2(t-s)}(X^{(k,l)}_s) \I(s^{-1/2} X^{(k,l)}_s \in S_2)] 
& \leq &  \frac{c_1 \log \log 2(t-s)}{\log 2(t-s)} P[ \tau_c>s] \\
& \leq & \frac{C\log \log t}{\log t} \p(s)
\end{eqnarray*}
which then contributes to the error term $\mathcal{E}_1$. 

The remaining case where $\{k,l\} = \{k,i\}$ yields similar error estimates, but seems to be  
somewhat fiddlier. 
We again condition at time $s$ to get 
\begin{equation} \label{hardercase}
\PR_[\tau_{c}>s, \tau_{ij}>t, |X^{(i,j)}_t|=1, \tau_{ki} \leq t]
= \E[\I(\tau_{c}>s) H_{t-s}(X^{(i)}_s,X^{(j)}_s,X^{(k)}_s)]
\end{equation}
where
\[
H_t(x_1,x_2,x_3) = \PR_{(x_1,x_2,x_3)}[ \tau_{12}>t, \tau_{13} \leq t, |X^{(1,2)}_t|=1].
\]
Our aim is to avoid an estimate that involves a three particle calculation, and we split $H$ into two parts
on each of which we will need only single particle estimates. 
We introduce another intermediate time $r = t \log^{-1}\!t$.
We split the expectation for $H_t = H^{(1)}_t + H^{(2)}_t$ according to $\tau_{1,3} \in [t-r,t]$ or $\tau_{1,3} \in [0,t-r]$. 
The first part is bounded as follows:
\begin{eqnarray*}
H^{(1)}_t(x_1,x_2,x_3) & : = & \PR_{(x_1,x_2,x_3)}[ \tau_{12}>t, \tau_{13} \in (t-r,t], |X^{(1,2)}_t|=1] \\
& \leq & \PR_{(x_1,x_2,x_3)}[ \tau_{13} \in (t-r,t], |X^{(1,2)}_t|=1] \\
& \leq & \frac{C}{t} \PR_{(x_1,x_3)}[ \tau_{13} \in (t-r,t]] \\
& = & \frac{C}{t} \left(G_{2t}(x_{13}) - G_{2(t-r)}(x_{13}) \right)
\end{eqnarray*}
where in the second inequality we have conditioned on the paths $X^{(1)}$ and $X^{(3)}$ and used a 
local central limit theorem to bound the resulting probability for $|X^{(2)}_t - z|=1$.
A generator calculation shows that
\begin{equation} \label{Gderivative}
\frac{d}{dt} G_{t}(x) = \PR_x [ \tau_0 > t, |X|_t =1] = F_t(x)
\end{equation}
so that using Lemma \ref{Gbounds} (b)
\[
H^{(1)}_t(x_1,x_2,x_3) \leq \frac{C}{t} \frac{r}{2(t-r) \log 2(t-r)} = O \left(\frac{1}{t \log^2 t} \right)
\]
and this part of (\ref{hardercase}) contributes to the error term $\mathcal{E}_1$. For $H^{(2)}$ we 
use Lemma \ref{Gbounds} (b) to find
\begin{eqnarray*}
&& \hspace{-.4in} H^{(2)}_t(x_1,x_2,x_3) \\
 & : = & \PR_{(x_1,x_2,x_3)} [ \tau_{12}>t, \tau_{13} \leq  t-r, |X^{(1,2)}_t|=1] \\
& = & \sum_{y=(y_1,y_2,y_3)} \PR_{(x_1,x_2,x_3)} [\tau_{1,2} > t-r, \tau_{13} < t-r, X_{t-r} = y]
\PR_{(y_1,y_2)} [\tau_{12} > r, |X^{(1,2)}_r|=1] \\
& \leq & \sum_{(y_1,y_2,y_3)} \PR_{(x_1,x_2,x_3)} [\tau_{13} < t-r, X_{t-r} = y]
\left( \frac{4c_2}{\log 2r} \PR_{(y_1,y_2)} [|X^{(1,2)}_r|=1]  + \frac{4c_2}{2r \log^2 2r} \right) \\
& = & \frac{4c_2}{\log 2r}  \PR_{(x_1,x_2,x_3)} [\tau_{13}<t-r, |X^{(1,2)}_t|=1] + 
\frac{4c_2}{2r \log^2 2r} \PR_{(x_1,x_2,x_3)} [\tau_{13}<t-r] \\
& \leq & \frac{C}{t \log t} \PR_{(x_1,x_3)} [\tau_{13} < t-r].
\end{eqnarray*}
We may again split the expectation (\ref{hardercase})  according to $s^{-1/2} X^{(k,i)}_s \in S_1 \cup S_2$ and when
$2\gamma +1 = \mu$, at the expense of a further 
contribution to $\mathcal{E}_2$, we need only consider $S_2$.
Then choosing $\beta > 2 \gamma + \alpha$ the estimate from Lemma \ref{Gbounds} (a) applies at large $t$ to show
\[
H^{(2)}_{t-s}(X^{(i)}_s,X^{(j)}_s,X^{(k)}_s) \I(s^{-1/2} X^{(k,i)}_s \in S_2) = O\left(\frac{1}{t \log^2 t}\right). 
\]
and the contribution from $H^{(2)}$ to (\ref{hardercase}) can be absorbed into $\mathcal{E}_1$.

\noindent
\textbf{Step 2.}  For the error in the second approximation  (\ref{app2}) we break the expectation 
\[
\E_x \left [ \I(\tau_c>s) \left(F_{2(t-s)}(X^{(i,j)}_s) - \frac{4}{2(t-s) \log 2(t-s)}\right) \right].
\]
into three parts according to whether $s^{-1/2} X^{(i,j)}_s \in S'_1 \cup S'_2 \cup S'_3$
where
\[
S'_1 = \{|z| \leq \log^{-\gamma_1}\!s \}, \quad
S'_2 = \{  \log^{-\gamma_1}\!s < |z| < \log^{\gamma_2}\!s \} , \quad
S'_3 =  \{ \log^{\gamma_2}\!s \leq |z|\} .
\]
The part $S'_1$ as before leads to a contribution to the error term $\mathcal{E}_2$ if we take
$2 \gamma_1 + 1=\mu$.
Large deviation estimates (see (\ref{LD})) show, when $\gamma_2 > 0$, the part $S'_3$ is $O(\log^{-\mu}\!t)$ and 
also contributes only to the error term $\mathcal{E}_2$. 
When $s^{-1/2} X^{(i,j)}_s \in S'_2 $ we have
\[
\frac{t}{(\log t)^{\alpha+2\gamma_1}}
 \leq s (\log s)^{-2 \gamma_1} \leq |X^{(i,j)}_s|^2  \leq  s (\log s)^{2 \gamma_2} 
\leq \frac{t}{(\log t)^{\alpha-2\gamma_2}},
\]
We choose $\alpha> 2 \gamma_2+1$. Then, by taking $\beta_1 > \alpha + 2 \gamma_1$ and $\beta_2 \in [1, \alpha - 2 \gamma_2)$, 
we may apply the estimates in Lemma \ref{Gbounds} (c) when  $t$ is large to see that
\begin{eqnarray*}
\left|F_{2(t-s)}(X^{(i,j)}_s) - \frac{4}{2(t-s) \log 2(t-s)}\right| \I(s^{-1/2} X^{(i,j)}_s \in S'_2)  & \leq & 
\frac{4c_4 \log \log 2(t-s)}{2(t-s) \log^2 \!2(t-s)} \\
& = & O \left( \frac{\log \log t}{t \log^2 \!t} \right)
\end{eqnarray*}
and we find another error term contributing to $\mathcal{E}_1$. 

\textbf{Step 3.} The final approximation (\ref{app3}) is straightforward using the definition of $s$ and,
since we take $\alpha \geq 1$,  
adjusts the error $\mathcal{E}_1$ without destroying the desired error bound.  This completes the proof
of the errorbounds (\ref{errorbounds}).

The large deviation estimates we use above can all be deduced from the following standard one dimensional estimate. 
One coordinate of one of our particles it is a rate $1/2$ simple random walk on $\Z$. Writing $(Z_t)$ for such a process 
started at $0$ we have for some $c_5,c_6>0$ 
\begin{equation} \label{LD}
\PR[\sup_{s \leq t} |Z_s| \geq t^{1/2} \log^r \!t] \leq c_5 \exp(- c_6 \log^{2r}\!t) \quad \mbox{for all $t>0$.}
\end{equation}
This can be established, as usual, at a fixed $t$ by Chebyshev's inequality using exponential moments, 
and a reflection principle can be used to control the supremum.
We use it when $r>1/2$ to obtain a bound that is $O(t^{-\mu})$ for any $\mu>0$, and when $r>0$ to obtain a 
bound that is $O(\log^{-\mu}\!t)$ for any $\mu>0$.

\textbf{Proof of Lemma \ref{Gbounds}.} Known estimates on the hitting time of the origin 
for a discrete time simple random walk on $\Z^2$
(which we denote by $\hat{\tau}_0$) include
\[
\PR_{\pm e_i}[ \hat{\tau}_0 \geq n]  =  \frac{\pi}{\log n} + O((\log n)^{-2})
\]
 and
\[
\PR_x[ \hat{\tau}_0 \leq n]  = \frac{\log(n/|x|^2)}{\log n} \left(1+ O\left(\frac{\log \log (n/|x|^2)}{\log (n/|x|^2)}\right)  \right) \qquad
 \mbox{when $n^{1/6} \leq |x| \leq \frac{1}{20} n^{1/2}$,}
\]
(see R\'ev\'esz \cite{revesz} Lemma 20.1 and Theorem 20.3).
It is straightforward to Poissonize the number of jumps to establish part (a) of the Lemma.

For part (b) we make use the simple bound $q_t(y,z) \leq p_t(z-y)$ and a local central limit theorem for 
$p_t$:  for any $\nu>3/2$
\begin{equation} \label{lclt}
p_t(x) = \frac{1}{\pi t} \exp(-|x|^2/t) + R(t,x) \quad \mbox{where $\;\;\sup_x R(t,x) = O(t^{-\nu})$.}
\end{equation}
(This can be deduced from the discrete time version in \cite{revesz} Lemma  17.6). 

Splitting at a time $s = t \log^{-\alpha}\!t$ for some $\alpha >5$ we find
\begin{eqnarray*}
q_t(\pm e_i,z) &=&\sum_y q_s(\pm e_i,y) q_{t-s}(y,z) \\
& \leq & \sum_y q_s(\pm e_i,y) p_{t-s}(z-y) \\
& = &  \sum_y q_s(\pm e_i,y) p_t(z) + \mathcal{E}_5 \\
& = &  \PR_{e_i}[\tau_0 >s] \, p_t(z) + \mathcal{E}_5.
\end{eqnarray*}
Combined with part (a), this leads to part (b) of the lemma provided we can show $\mathcal{E}_5 = O(1/t \log^2 \!t)$.
We may restrict to summing over $|y| < s^{1/2} \log^r \!s$ for some $r>1/2$ by the large deviation estimate (\ref{LD}). 
Then (\ref{lclt}) and simple estimates on the derivatives of the Gaussian transition density show
\begin{eqnarray*}
|p_{t-s}(z-y) - p_t(z)| & \leq & |p_{t-s}(z-y) - p_t(z-y)| + |p_{t}(z-y) - p_t(z)| \\
& \leq &  C s (t-s)^{-2} + C |y| t^{-3/2} + O(t^{-\nu})
\end{eqnarray*}
which leads to the desired bound.  

Repeating the last argument but tracking the upper bound more carefully (using (\ref{lclt}) and part (a)) we have
\begin{eqnarray*}
\sup_z q_t(\pm e_i,z)  & \leq & \|p_{t-s}\|_{\infty} \PR_{\pm e_i} [\tau_0>s] \\
& = & \left( \frac{1}{\pi(t-s)} + O(t^{-\nu})\right) \left(\frac{\pi}{\log s} + O\left(\frac{1}{\log^{2}\!s}\right) \right) \\
& = &  \frac{1}{t \log t} \left( 1 + O\left(\frac{1}{\log t}\right) \right).
\end{eqnarray*}
This proves the upper bound for part (c). 
For the lower bound we use
\begin{equation} \label{splitlem}
q_t(\pm e_i,z)  =  \PR_{\pm e_i}[\tau_0 > t, X_t =z] =  \PR_{\pm e_i}[\tau_0 >s, X_t =z] - \PR_{\pm e_i}[s<\tau_0 <t, X_t =z] 
\end{equation}
The first term in (\ref{splitlem}) is
\[
\PR_{\pm e_i}[\tau_0 >s, X_t =z] = \sum_y q_s(\pm e_i,y) p_{t-s}(z-y) 
\]
and we may again restrict to $|y| < s^{1/2} \log^r \!s $ for some $r>1/2$ by a large deviation tail estimate.
Recall that we are considering $ |z|^2 \leq \frac{t}{\log^{\beta_2}\!t}$ where $\beta_2 \geq 1$. 
The local central limit theorem (\ref{lclt}) can be used, when $\alpha \geq 2r+1$, to show for these $y,z$ that $p_{t-s}(y,z) = \frac{1}{\pi t} + O(1/t \log t)$
and then
\[
\PR_{\pm e_i}[\tau_0 >s, X_t =z] = \PR_{\pm e_i} [\tau_0>s] \left(\frac{1}{\pi t} + O\left(\frac{1}{t \log t}\right)\right) = \frac{1}{t \log t} + O\left(\frac{1}{t\log^2\!t} \right).
\]
It remains to bound the second term in (\ref{splitlem}):
\begin{equation} \label{finalerror}
\PR_{\pm e_i}[s <\tau_0<t, X_t =z] = \E_{\pm e_i}[G_{t-s}(X_s,z) \I(s <\tau_0)]
\end{equation} 
where $G_t(y,z) = \PR_y[\tau_0 <t, X_t =z]$. We may add the indicator
$\I(|X_s| \geq s^{1/2} \log^{-\gamma}\!s)$ into (\ref{finalerror}) for
$\gamma \geq 1$ at the expense only of an $O(1/t \log^2 \!t)$ error. Then we use
\begin{eqnarray*}
G_{t}(y,z) &=&  \PR_y[\tau_0 \in [0,t/2), X_t =z] + \PR_y[\tau_0 \in [t/2,t], X_t =z] \\
& \leq & \PR_y[\tau_0 \in [0,t/2], X_t =z] + \PR_z[\tau_0 \in [0,t/2], X_t =y] \\
& \leq & \frac{C}{t} \left( \PR_y[\tau_0 \leq t/2] + \PR_z[\tau_0 \leq t/2] \right).
\end{eqnarray*}
Using this bound, part (a) of the Lemma shows (taking $\beta > \max\{\beta_1, 2 \gamma + \alpha\}$) that 
\[
|G_{t-s}(y,z)| = O\left(\frac{\log \log t}{t \log t}\right) \qquad \mbox{for
$|z|^2 \geq t \log^{-\beta_1}\!t$ and $|y|^2 \geq s \log^{-2\gamma}\!s$}
\]
and this is sufficient to show that
(\ref{finalerror}) is $O(\log \log t/t \log^2 \!t)$, completing the proof of the lower bound. 

Under our restrictions that $|z|^2 \geq t \log^{-\beta_1}\!t$ and $|y|^2 \geq s \log^{-2\gamma}\!s$
part (a) of the Lemma shows (taking $\beta > \max\{\beta_1, 2 \gamma + \alpha\}$) that 
$|G_{t-s}(y,z)| = O(\log \log t/t \log t)$.  This is sufficient to show that
(\ref{finalerror}) is $O(\log \log t/t \log^2 \!t)$ which completes the proof of the lower bound. 
\subsection{Particle intensity $\rho_1$ for infinite particle system} \label{s2.2}
We may construct the particle system as an infinite system of SDEs driven by Poisson processes. Let $(P(x,y): x,y \in \Z^2)$ be an I.I.D. family of rate $\frac14$ Poisson processes.
Write $x \sim y$ when $x$ and $y$ are neighbours, that is when $|x-y|=1$. The process $t \to P_t(x,y)$ will trigger the jumps of particles
from $x$ to $y$. The process $(\xi_t(x): x \in \Z^2)_{t \geq 0}$ is the unique process with values in $\{0,1\}^{\Z^2}$ solving for all $x$ the system
\begin{equation} 
d \xi_t(x) =  - \sum_{y \sim x} \xi_{t-}(x) dP_t(x,y) + \sum_{y \sim x} (1-\xi_{t-}(x)) \xi_{t-}(y) dP_t(y,x) 
 \quad \mbox{for all $x \in \Z^2,\,t \geq 0$} \label{system}
\end{equation} 
with initial condition $\xi_0 \equiv 1$.  
This is very close to a graphical construction, but the system can also be treated in a standard 
differential equations manner. 

The particle intensity $\rho_1(t) =  \E[\xi_t(x)] $ is independent of $x$ due to the choice of the initial conditions. 
When we consider $N$ positions we consider disjoint sites and write 
\[
x= (x_1,\ldots,x_N) \in \Z^{2N}_d := \{z \in (\Z^2)^N: z_i \neq z_j \mbox{ for all $i \neq j$}\}.
\]
Then
\[
\rho_N(x;t) = \PR[\mbox{$x_1,\ldots,x_N$ all occupied at time $t$}] = \E[\xi_t(x_1) \ldots \xi_t(x_N)]
\quad \mbox{for $x \in \Z^{2N}_d$.} 
\]
The variables $(\xi_t(x):x \in \Z^2)$ have negative correlation properties, in particular (see Lemma 2.5 in \cite{KvdB2}) 
\begin{equation} \label{NC}
\rho_N(x;t) \leq \rho_1^N(t) \quad \mbox{for all $t \geq 0,x \in \Z_d^{2N}$ and $N \geq 2$.}
\end{equation}
We will also use
\begin{equation} \label{NCplus}
\E[ \xi_t(x_1) \xi_t(x_2) \xi_t(x_3)] \leq \E[ \xi_t(x_1)] E[\xi_t(x_2) \xi_t(x_3)] 
\end{equation}
which can be established in the same manner as Lemmas 2.4-2.9
in \cite{KvdB2}, where they are deduced from the BKR inequality for an approximating discrete time process. 

\noindent
\textbf{Notation.}
For $f: \Z^2 \to \R$ and $\psi: (\Z^2)^2 \to \R$ we write, when these are well defined, 
\[
\langle \xi_t,f \rangle = \sum_x \xi_t(x) f(x) 
\qquad \langle \xi_t * \xi_t, \psi \rangle = \sum_{x,y} \xi_t(x) \xi_t(y) \psi(x,y).
\]
We use the two particle test function 
\begin{equation} \label{testfn}
\psi_t(x,y) = \PR_{(x,y)}[\tau_c >t, X^{(1)}_t=0, X^{(2)}_t = e_1]. 
\end{equation}
Recall we are writing $p_t(y)$ for the transition probabilities for a simple random walk on $\Z^2$. 
The effective rate equation for $\rho_1$ is found by the following approximations: 
\begin{eqnarray}
\frac{d}{dt} \rho_1(t) & = & - \E [ \xi_t(0) \xi_t(e_1)] \label{B1} \\
& = & - \E [ \langle \xi_{t-s} * \xi_{t-s}, \psi_s] + \mathcal{E}_1(t) \label{B2} \\
& = & - \frac{\pi}{\log s} \E [ \langle \xi_{t-s}, p_s \rangle^2]  + \mathcal{E}_2(t)  \label{B3} \\
& = & - \frac{\pi}{\log s} \left( \E [ \langle \xi_{t-s}, p_s \rangle] \right)^2 
 + \mathcal{E}_3(t)  \label{B4} \\
& = & - \frac{\pi}{\log s} \rho_1^2 (t-s) 
 + \mathcal{E}_3(t)  \nonumber \\
& = &  - \frac{\pi}{\log s} \rho_1^2(t) + \mathcal{E}_4(t). \label{B5}
\end{eqnarray}
The error bounds below will show that 
\begin{equation} \label{ERROR}
|\mathcal{E}_4(t)|  \leq C \left( \frac{\rho_1(t-s)}{s \log s} + \frac{\rho_1^2(t-s)}{\log^{\frac32}\!s} +\frac{s \rho_1^3(t-2s)}{\log^{2}\!s} \right) 
\quad \mbox{for all $2 \leq  s  \leq t/2$.} 
\end{equation}
Combined with the known crude upper and lower estimates for $\rho_1(t)$ 
(see \cite{Bramson+G} equations (14), (15) and (20))
\begin{equation} \label{crude}
0 < c_7 \frac{\log t}{t} \leq \rho_1(t) \leq c_8 \frac{\log t}{t} \quad \mbox{for all $t \geq 2$}
\end{equation}
and the choice $s = t \log^{-\frac12} \!t$ we find 
\[
\frac{d}{dt} \rho_1(t) =   - \frac{\pi}{\log t} \rho_1^2(t) + O\left(\frac{\log^{\frac12}\!t}{t^2}\right).
\]
It is then a calculus exercise (left for the reader) to derive the desired asymptotics (\ref{sawyer}) from this and the crude bounds (\ref{crude}).

We will now bound each of the increments in the error term 
$D \mathcal{E}_k = \mathcal{E}_k - \mathcal{E}_{k-1}$ for $k=1,2,3,4$. 
We need some detailed information about the two particle function test function $\psi$, namely the following analogue of
Lemma 12 in \cite{KvdB1}: for any $\epsilon>0$ 
\begin{equation} \label{Jamielemma}
\left| \sum_{x,y} \psi_t(x,y)  - \frac{\pi p_t(x) p_t(y)}{\log t}  \right|  = O\left( \frac{1}{\log^{2-\epsilon}\!t}\right) . 
\end{equation}
This is a special case of Lemma \ref{Jamieplus} below. We split the error in (\ref{B3}) into two:
\[
D\mathcal{E}_2(t)  = \sum_{x \neq y}  \E[ \xi_{t-s}(x) \xi_{t-s}(y)] \left( \frac{\pi p_s(x) p_s(y) }{\log s}  -\psi_s(x,y) \right) \\
+ \sum_{x}  \E[ \xi_{t-s}(x)]  \frac{\pi}{\log s} p_s^2(x). 
\]
Then using $\sum_x p^2_t(x) \leq \|p_t\|_{\infty} \leq C t^{-1}$, the bound (\ref{Jamielemma})  
and the negative correlation (\ref{NC}), we have
\[
|D\mathcal{E}_{2}(t) | \leq  
\frac{C}{(\log s)^{\frac32}} \rho_1^2(t-s)  + \frac{C}{s \log s} \rho_1(t-s) \quad \mbox{for all $2 \leq s \leq t$,}
\]
(we may choose $\epsilon = 1/2$ in (\ref{Jamielemma}) to match other error terms, as this term is not a biting error).  
The error bound in the variance estimate (\ref{B4}) is immediate in our setting via negative correlation since
\[
\mbox{Var}( \langle \xi_t, f \rangle ) = \sum_{x,y} f(x) f(y) \mbox{Cov}( \xi_t(x), \xi_t(y)) \leq \sum_x f^2(x) \mbox{Var}(\xi_t(x)) \leq 
 \sum_x f^2(x) \E[\xi_t(x)]
\]
using $\xi_t(x) \in \{0,1\}$ in the last step. This implies the error bound 
\[
|D\mathcal{E}_3(t)| \leq  \frac{\pi}{\log s} \rho_1(t-s) \langle p_s^2,1 \rangle  \leq C  \frac{\pi}{s \log s}\rho_1(t-s)
\quad \mbox{for all $2 \leq s \leq t$.}
\]
Compensating the Poisson processes, that is replacing $dP_t(x,y)$ by $(dP_t(x,y) - \frac14 dt) + \frac14 dt$, we collect the 
martingale increments into a single martingale $M_t(x)$  to find
\[
d \xi_t(x) = \Delta \xi_t(x) dt - \xi_t(x) S\xi_t(x)  dt + dM_t(x)
\]  
where 
\[
S f(x) = \frac14 \sum_{y \sim x} f(y), \qquad \Delta f(x) = Sf(x) - f(x)
\]
so that $\Delta$ is the standard discrete Laplacian. 
The initial conditions imply that $\E[\Delta \xi_t(x)] = 0$ and $\E[\xi_t(x) \xi_t(x \pm e_i)] = \E[ \xi_t(0) \xi_t(e_1)]$ so that taking expectations we  reach (\ref{B1}). 

A longer, but similar, exact calculation shows that 
\begin{equation}
d(\xi_t(x) \xi_t(y))  =  \xi_t(x) \Delta \xi_t(y) dt + \xi_t(y) \Delta \xi_t(x) dt - \mu^{(1)}_t(x,y) dt + dM_t(x,y)  \label{2pta}
\end{equation}
for a martingale $M_t(x,y)$, where 
\begin{eqnarray*}
\mu^{(1)}_t(x,y) & = & \xi_t(x) \xi_t(y) (S\xi_t(x) + S\xi_t(y))  \nonumber \\
&& \hspace{.2in} + \frac14 \I(x \sim y) \left( \xi_t(x) + \xi_t(y) - 2 \xi_t(x) \xi_t(y) \right) \nonumber \\
&& \hspace{.4in} - \I(x=y) \left( \xi_t(x) + (1-\xi_t(x)) S\xi_t(x) \right).
\end{eqnarray*}
Suppose $\phi: [0,t] \times (\Z^2)^2 \to \R$ satisfies $ \sup_{s \leq t} \sum_{x,y} |\phi_s(x,y)| + | \dot{\phi}_s(x,y)| < \infty$. Then 
developing $d(\phi_t(x,y) \xi_t(x) \xi_t(y))$, summing over $x,y$, and applying discrete
integration by parts gives
\begin{equation}
d \langle \xi_t * \xi_t, \phi_t \rangle  =  \langle \xi_t * \xi_t, \Delta_{x,y} \phi_t + \dot{\phi}_t \rangle dt 
- \langle \mu^{(1)}_t, \phi_t \rangle dt + \langle dM_t,\phi_t \rangle. 
 \label{2ptb}
\end{equation}
Here $\Delta_{x,y}$ means $\Delta$ applied in both $x$ and $y$ variables. The integrability of $\phi$ implies
that the sum $ \langle dM_t,\phi_t \rangle$ is still a martingale.  

We apply this to the function $\phi_s(x,y) = \psi_{t-s}(x,y)$ from (\ref{testfn}). Note that 
\[
\dot{\psi}_s = \Delta_{x,y} \psi_s \quad \mbox{for $s \geq 0$ and $x \neq y$}
\]
while $\psi_s(x,x) = 0$ for all $s \geq 0$. For this choice (\ref{2ptb}) becomes 
(after combining the term $\Delta \psi_s(x,x)$ carefully) exactly the term 
\begin{equation}
d \langle \xi_s * \xi_s, \psi_{t-s} \rangle  = - \langle \mu^{(2)}_s, \psi_{t-s} \rangle dt + \langle dM_s,\psi_{t-s} \rangle
 \label{2ptc}
\end{equation}
where
\[
\mu^{(2)}_t (x,y) = \left\{ \begin{array}{ll}
\xi_t(x) \xi_t(y) (S\xi_t(x) + S\xi_t(y)) & \mbox{if $|x-y|>1$,} \\
\frac14 \xi_t(x) \xi_t(y)  \sum_{w \neq x,y} \xi_t(w) \I(\mbox{$w \sim x$ or $w \sim y$})
& \mbox{if $x \sim y$.}
\end{array} \right.
\]
Note the initial condition $\psi_0(x,y) = \I(x=0,y=e_1)$. Since $\mu^{(2)}_t(x,y) \geq 0$ for $x \neq y$
and $\psi_s(x,x)=0$, integrating (\ref{2ptc}) over $[t-s,t]$ and taking expectations yields
\[
\E[\xi_t(0) \xi_t(e_1)] \leq  \E[\langle \xi_{t-s} * \xi_{t-s}, \psi_{s} \rangle] \leq  \langle 1, \psi_s \rangle \rho_1^2(t-s).
\]
Then the bound (\ref{Jamielemma}) implies (the initial conditions imply this is the same for any pair $x \sim y$)
\begin{equation}  \label{newestimate}
\E[\xi_t(x) \xi_t(y)] \leq  \frac{C}{\log s} \rho_1^2(t-s) \quad \mbox{for all $2 \leq s \leq t$ and $x \sim y$.} 
\end{equation}
We remark that using just negative correlation at time $t$ yields $\E[\xi_t(x) \xi_t(y)] \leq \rho_1^2(t)$, and the improvement given 
by (\ref{newestimate}) is crucial below. 

Furthermore, the formula for $\mu^{(2)}_t(x,y)$ involves $\xi_t(x)\xi_t(y)\xi_t(w)$ where $w$ is disjoint from $\{x,y\}$ and 
a neighbour of one of them. For example if $w \sim y$ we can bound by negative association (\ref{NCplus})
and then employ (\ref{newestimate}). This implies that we can bound
\begin{equation} \label{mubound}
\E[ \hat{\mu}_t(x,y)] \leq \frac{C}{\log s} \rho_1(t) \rho_1^2 (t-s)  \leq \frac{C}{\log s} \rho_1^3(t-s) \quad 
\mbox{for all $2 \leq s \leq t$ and $x \neq y$.} 
\end{equation}
Using Lemma \ref{Jamielemma} we have $\langle \psi_r,1 \rangle \leq C \min\{1, \log^{-1} \!r\}$, so 
we can estimate the error in the first approximation step (\ref{B2}): 
 \begin{eqnarray*}
|\mathcal{E}^{(1)}_t| & \leq & \int^t_{t-s} \langle E[\mu^{(2)}_{r}], \psi_{t-r} \rangle dr \\
& \leq &   \frac{C}{\log s} \rho_1^3(t-2s)
\int^t_{t-s} \langle 1, \psi_{t-r} \rangle dr \\
& \leq &   \frac{Cs}{(\log s)^2} \rho_1^3(t-2s)
\quad \mbox{for all $2 \leq s \leq t/2$.} 
\end{eqnarray*}
The final approximation step (\ref{B5}) replaces $\rho_1^2(t-s)$ by $\rho_1^2(t)$, which by (\ref{newestimate}) creates an error
\begin{eqnarray*}
\left| \rho_1^2(t-s) - \rho_1^2(t) \right| & \leq & 
2 \rho_1(t-s) (\rho_1(t-s) - \rho_1(t)) \\ 
& = & 2 \rho_1(t-s) \int^t_{t-s} E[\xi_r(0) \xi_r(e_1)] dr \\
& \leq & C \frac{s}{\log s}  \rho_1(t-s) \rho^2_1(t-2s) \quad \mbox{for all $2 \leq s \leq t/2$.} 
\end{eqnarray*}
This leads to an error increase $D\mathcal{E}_4$ of the same order as $\mathcal{E}_1$.
Summing all the error terms leads to  (\ref{ERROR}). 

\subsection{Higher correlations $\rho_N$ for infinite particle system} \label{s2.3}

We need a graphical construction of the system that allows us to track the trajectories of the particles and identify when particles 
at different times are descended from each other. 
Indeed the process can be constructed inductively by adding successive random walk paths in any order, 
coalescing any path that meets a previously defined path. This allows us to write, for $x,y \in \Z^{2}$ and
$0\leq s<t$
\[
(y,s) \to (x,t)
\]
to mean that the particle at positions $y$ at time $s$ moves to position
$x$ at time $t$. 

Fix $0< t-s< t$ and $x = (x_1,\ldots,x_N) \in \Z^{2N}_d$. Define for $y=(y_1,\ldots,y_N)$ 
\[
\Omega_y = \{ \xi_{t-s}(y_i) =1, (y_i,t-s) \to (x_i,t) \; \mbox{for $i=1,\ldots,N$}\}.
\]
All we require from the graphical construction is that $\rho_N(x;t) = \PR \left[ \bigcup_{y \in \Z^{2N}_d} \Omega_y \right]$.
We define an $N$ particle test function
\[
\psi^{(N)}_t(x,y) = \PR_{y}[X_t = x,\tau_c>t]
\]
and note that $\psi^{(N)}_t(x,y) = \psi^{(N)}_t(y,x)$. 
Using the Markov property at time $t-s$ and the 
negative correlation (\ref{NC}) we have 
\begin{eqnarray*}
\rho_N(x;t) & = & \PR \left[ \bigcup_{y \in \Z^{2N}_d} \Omega_y \right] \\
& \leq & \sum_{y \in \Z^{2N}_d}
 \PR[ \Omega_y] \\
& = &\sum_{y \in \Z^{2N}_d} \rho_N(y;t-s)  \psi^{(N)}_{s}(y,x) \\
& \leq &  \rho_1^N(t-s) \sum_{y}   \psi^{(N)}_{s}(y,x) \\
& = &  \rho_1^N(t-s) \sum_{y}  \psi^{(N)}_{s}(x,y) \\
& = &  \rho_1^N(t-s) \PR_{x} [\tau_c  >s]. 
\end{eqnarray*}
We choose $s = t \log^{-\frac12}\! t$ and then 
the asymptotics (\ref{sawyer}) for $\rho_1$ and (\ref{MR1})  for $\p(s)$ give
\[
\rho_N(x;t) \leq c_0(x)  (\log t)^{-{N \choose 2}} \rho_1^N(t)  \left(1+ O\left(\frac{1}{\log^{\frac12}\!t}\right) \right).
\]
To show this upper bound is correct we have, informally, to complete two remaining steps: 
(\textbf{Step 1}) show our estimate on $ \sum_y P[ \Omega_y] $ is accurate; (\textbf{Step 2}) use inclusion 
exclusion to show we have not lost anything in the union bound. These two steps will show the lower bound
and establish the multiplicative error term of the order $(1+O(\log^{-(\frac12-\delta)}\!t))$. Throughout
we are taking $s = t \log^{-\frac12}\! t$ and consider only $t$ so that $1\leq s \leq t/2$. 

\noindent
\textbf{Step 1.} We need the following estimates for the $N$ particle test function $\psi^{(N)}$, 
which are proved at the end of this section. Recall for $y \in \Z^{2N}$ that $p_t(y) = \prod_{i=1}^N p_t(y_i)$ is
the transition density for $N$ independent simple random walks. 
\begin{Lemma} \label{Jamieplus}
Fix $x$. 
For any $\epsilon>0$ 
\begin{equation} \label{JLplus}
\sum_{y} \left|\psi^{(N)}_t(x,y) - c_0(x) (\log t)^{-{N \choose 2}} p_t(y) \right| = O \left(\frac{1}{(\log t)^{{N \choose 2}+1-\epsilon}} \right). 
\end{equation}
Moreover, for any $\mu>0$, 
\begin{equation} \label{UBNplus}
\psi^{(N)}_t(x,y) \leq c_9(x) (\log t)^{-{N \choose 2}} p_t(y) + c_{10}(x,\mu) t^{-N} \log^{-\mu}\!t
\quad \mbox{for all $t \geq 2$.}
\end{equation}
\end{Lemma}
Using (\ref{JLplus})  we may write
\begin{eqnarray*}
 \sum_{y \in \Z^{2N}_d} \PR[ \Omega_y] & = & \sum_{y \in \Z^{2N}_d} \rho_N(y;t-s) \psi^{(N)}_{s}(y,x) \\
 & =  & c_0(x) (\log s)^{-{N \choose 2}} \sum_{y \in \Z^{2N}_d} \rho_N(y;t-s) p_s(y)
 +O \left(\frac{\rho_1^N(t)}{(\log t)^{{N \choose 2}+1-\epsilon}}\right).
 \end{eqnarray*}
 We rewrite 
 \begin{eqnarray*}
&& \hspace{-.5in}  \sum_{y \in \Z^{2N}_d} \rho_N(y;t-s) p_s(y)  \\
& = & \sum_{y \in \Z^{2N}} \rho_N(y;t-s) p_s(y)- \sum_{y \in \Z^{2N} \setminus \Z^{2N}_d} \rho_N(y;t-s) p_s(y) \\
 &  = & E[ \langle \xi_{t-s}, p_s \rangle^N]  
 - \sum_{y \in \Z^{2N} \setminus \Z^{2N}_d} \rho_N(y;t-s) p_s(y)  \\
 & \geq &  \left(E[ \langle \xi_{t-s}, p_s \rangle]\right)^N  
  - \sum_{y \in \Z^{2N} \setminus \Z^{2N}_d} \rho_N(y;t-s) p_s(y) \\
  & = & \rho_1^N(t-s) 
   - \sum_{y \in \Z^{2N} \setminus \Z^{2N}_d} \rho_N(y;t-s) p_s(y)
\end{eqnarray*}
using Jensen's inequality in the penultimate step.
We must now show that the sum $\sum_{y \in \Z^{2N} \setminus \Z^{2N}_d}$ gives a subleading contribution. This sum
can be organised according to the number $K<N$ of distinct coordinates in $y$. The terms with $K$
distinct coordinates give contributions of the form 
\[
 \sum_{y \in \Z^{2K}_d} \E[\xi_{t-s}(y_1)\ldots \xi_{t-s}(y_K)] \prod_{i=1}^K p^{n_i}_s(y_i) 
\]
where $1 \leq n_i \in \N$ satisfy $\sum n_i = N$. Using $\langle p^n_s,1 \rangle \leq C s^{1-n}$ this 
is bounded by
\[
C  \rho_1^K(t-s) s^{K - N} 
\leq C  \frac{\log^K\!t}{t^K} \frac{t^{K - N}} {(\log t)^{(K-N)/2}}
= C   \frac{\log^N\!t}{t^N} (\log t)^{-(N-K)/2}.
\]
which is small compared to the main term $\rho_1^N(t)$. Collecting the arguments in this step we find
\[
 \sum_{y \in \Z^{2N}_d} \PR[ \Omega_y] \geq 
 c_0(x)  (\log t)^{-{N \choose 2}} \rho_1^N(t)  \left(1- O\left(\frac{1}{\log^{\frac12}\!t}\right) \right).
\]

\noindent
\textbf{Step 2.} We apply the Bonferroni inequality
\[
\rho_N(x;t) \geq \sum_{y \in \Z^{2N}_d} \PR[\Omega_y]
- \frac12 \sum_{y \neq y'\in \Z^{2N}_d} \PR[\Omega_y \cap \Omega_{y'}]
\]
and aim to show that the second term is asymptotically subleading - that is it will be 
$O(t^{-N} (\log t)^{N - {N \choose 2} - \frac12+\delta})$.  

We break the double sum according to the number $K= \sum_{i \leq N} \I(y_i \neq y'_i)$ of distinct coordinates between
$y$ and $y'$.  Since $y \neq y'$ we have $1 \leq K \leq N$. It is enough to show that
sums of the form
\begin{equation} \label{lastsum}
\rho_1^{N+K}(t-s) \sum_{y \in \Z^{2N}_d} \sum_{y' \in \Z^{2K}_d} \PR[\OH_{y,y'}] I(y_i \neq y'_i, i = 1,\ldots,K)
\end{equation}
are asymptotically smaller, where
\[
\OH_{y,y'} = \{ (y_i,0) \to (x_i,s) \; \mbox{for $i =1,\ldots,N$}\}  \cap \;  \{ (y'_i,0) \to (x_i,s) \; \mbox{for $i =1,\ldots,K$}\}.
\]
Indeed, by varying $K \in \{1,\ldots,N\}$ and replacing $(x_1,\ldots,x_K)$ by other subsets
of $(x_1,\ldots,x_N)$ we find a finite number of events, whose sum dominates
$\sum_{y \neq y'} \PR[\Omega_y \cap \Omega_{y'}]$.

To estimate $\PR[\OH_{y,y'}]$ we may build a coalescing system from $N+K$
independent simple random walks $(X_r,X'_r)=((X^{(1)}_r,\ldots,X^{(N)}_r),(X^{'(1)}_r,\ldots,X^{'(K)}_r))$ 
over $r \in [0,s]$, started at $(y,y')=((y_1,\ldots,y_N),(y_1,\ldots,y_K))$ under $\PR_{y,y'}$. 
We may build a coalescing system by erasing paths after collisions (the details will not be useful for us). 
Let $\sigma_j$  be the collision time between $X^{(j)}$ and $X^{'(j)}$.
On the event  $\OH_{y,y'}$ we have $\sigma_j \leq s$ for all $j \leq K$. Furthermore by 
ignore other collisions we obtain an upper bound in terms of one and two particle probabilities:
\begin{eqnarray}
\PR[\OH_{y,y'}] & \leq & \PR_{y,y'}[X^{(i)}_s = x_i, \; \sigma_j \leq s, \, \mbox{for $i \leq N, j \leq K$}] \nonumber \\
& = & \prod_{i=1}^K \PR_{(y_i,y'_i)}[ X^{(1)}_s=x_i, \tau_{12} \leq s] \prod_{j=K+1}^N  \PR_{y_i}[ X^{(1)}_s=x_i]
. \label{pairbound}
\end{eqnarray}
The broad aim is to show that each of the $K$ collision events contributes $O(1/\log t)$ - at least for typical
configurations where $y_i,y'_i$ are of $O(s^{1/2})$ and are not too close. 

We start by making three reductions, shaving off various untypical parts of (\ref{lastsum}).
The first reduction is that we claim we may, for any $\epsilon>0$, restrict $y'$ in the double sum (\ref{lastsum}) to lie the box
\begin{equation} \label{Bbox}
B^{(K)}_s = \{z \in \Z^{2K}_d: \max_{i} |z_i| \leq s^{1/2} \log^{\epsilon}\!s\}.
\end{equation}
The large deviation estimate (\ref{LD}) shows that a particle starting outside $B^{(K)}_s$ has 
$O(\log^{-\mu}\!s)$ probability of reaching $x_k$ in time $s$. We claim that this will lead to the estimates,
for any $x \in \{x_1,\ldots,x_N\}$,
\begin{eqnarray} 
&& \sum_{y \in \Z^2} \sum_{y' \in \Z^2} \PR_{(y,y')}[ X^{(1)}_s=x, \tau_{12} \leq s] 
\leq 2s \label{claim2000a} \\
&& \sum_{y \in \Z^2} \sum_{y' \in \Z^2} \PR_{(y,y')}[ X^{(1)}_s=x, \tau_{12} \leq s] 
\I(|y'| > s^{1/2} \log^{\epsilon}\!s) = O(s \log^{-\mu}\!s) \label{claim2000b}
\end{eqnarray}
for any $\mu>0$. We prove this claim at the end of step 2. 
For $y' \not \in B_s^{(K)}$ at least one coordinate must satisfy $|y'_i| >s^{1/2} \log^{\epsilon}\! s$. By 
using the claims (\ref{claim2000a},\ref{claim2000b}), on the terms in (\ref{pairbound}) we find that
\[
\rho_1^{N+K}(t-s)  \sum_{y \in \Z^{2N}_d} \sum_{y' \in \Z^{2K}_d \setminus B^{(K)}_s} \PR[\OH_{y,y'}]
\leq C(\mu)  \rho_1^{N+K}(t-s) s^{K} \log^{-\mu}\!s.
\]
By taking $\mu$ large enough this will not affect the asymptotics (the leading term or the leading error term) for $\sum_y \PR[\Omega_y]$.
A very similar argument shows we may also reduce the sum over $y$ in (\ref{lastsum}) to $y \in B^{(N)}_s$. 

The second reduction is to restrict (\ref{lastsum}) to a further subset of $y'$:
\[
A_s = \{z \in \Z^{2K}_d: \min_{i} |z_i-y_i| \geq s^{1/2} \log^{-\gamma}\!s\}.
\]
When $y' \in A_s$ all the distances $|y'_i-y_i|$ are large enough that coalescence 
between the particles started at $y_i,y'_i$ is (logarithmically) unlikely. 
We can crudely bound
\begin{eqnarray*}
 \sum_{y \in B^{(N)}_s} \sum_{y' \in B^{(K)}_s \setminus A_s} \PR[\OH_{y,y'}] 
 & \leq &   \sum_{y \in \Z^{2N}} \sum_{y' \in B^{(K)}_s \setminus A_s} \PR_y[X_s=x]  \\
& \leq  & |B^{(K)}_s \setminus A_s|. 
\end{eqnarray*}
For $y' \in  B_s \setminus A_s$ at least one co-ordinate must lie in a box of width $2s^{1/2} \log^{-\gamma}\!s$, 
while the others may lie in $B_s$. This bounds the cardinality $ |B_s \setminus A_s|$ by 
$C s^K (\log s)^{2(K-1)\epsilon - 2 \gamma} $ and by taking $\gamma$ large this term also does not affect the asymptotics.

The third reduction is to remove the part of (\ref{lastsum}) where any of the collision times $\sigma_j$ occurs late in the interval $[0,s]$.
Let $r = t \log^{-\eta}\!t$, where $\eta>1/2$ so that $r < s$. We consider the sum
\begin{eqnarray}
&& \hspace{-.5in} \sum_{y \in B^{(N)}_s} \sum_{y' \in B^{(K)}_s}  \PR[\OH_{y,y'} \cap \{\sigma_1 > s-r\}] \label{sig1} \\
& \leq &\sum_{y \in B^{(N)}_s} \sum_{y' \in B^{(K)}_s}
 P_{(y_1,y'_1)} [X^{(1)}_s =x_1, \tau_{12} \in (s-r,s]] \prod_{i=2}^N p_s(x_i-y_i) \nonumber \\ 
 & \leq & C (s \log^{2\epsilon}\!s)^{K-1}  \sum_{y_1 \in B^{(1)}_s} \sum_{y'_1 \in B^{(1)}_s}  
 P_{(y_1,y'_1)} [X^{(1)}_s =x_1, \tau_{12} \in (s-r,s]]. \nonumber
\end{eqnarray}
We argue as in (\ref{Gderivative}) to get the bound
\begin{eqnarray*}
 && \hspace{-.5in} \sum_{y \in B^{(1)}_s} \sum_{y' \in B^{(1)}_s} P_{(y,y')}[ X^{(1)}_s = x, \tau_{12} \in (s-r,s]] \\
&= & \sum_{y \in B^{(1)}_s} \sum_{y^{''} \in B^{(1)}_s-y} P_{(0,y^{''})}[ X^{(1)}_s = x-y, \tau_{12} \in (s-r,s]] \\
& \leq & \sum_{y \in \Z^2} \sum_{y^{''} \in B^{(1)}_{s}-B^{(1)}_s} P_{(0,y^{''})}[ X^{(1)}_s = x-y, \tau_{12} \in (s-r,s]] \\
& = & \sum_{y^{''} \in B^{(1)}_{s}-B^{(1)}_s} P_{(0,y^{''})}[\tau_{12} \in (s-r,s]]  \\
& = & \sum_{y^{''} \in B^{(1)}_{s}-B^{(1)}_s}\left(G_{2s}(y^{''}) - G_{2(s-r)}(y^{''}) \right) \\
& \leq & C \frac{r |B^{(1)}_{s}-B^{(1)}_s|}{2(s-r) \log 2(s-r)} \\
& \leq & C t (\log t)^{-\eta -1 + 2 \epsilon}.
\end{eqnarray*}
Together these bounds imply that, by taking $\eta$ large, the term (\ref{sig1}) also does not affect the asymptotics.
The same applies replacing $\sigma_1$ by $\sigma_j$ for any $j \leq K$. 

We have reduced to the main part of the sum (\ref{lastsum}), namely
\[
\sum_{y \in B^{(N)}_s} \sum_{y' \in B^{(K)}_s \cap A_s} \PR[\OH_{y,y'} \I(\max_i \sigma_i \leq s-r)].
\]
Using the Markov property to condition at the positions of the particles at time $s-r$, and (\ref{UBNplus}) from
Lemma \ref{Jamieplus}, we rewrite this as
\begin{eqnarray*}
&&  \hspace{-.3in}
\sum_{y \in B^{(N)}_s} \sum_{y' \in B^{(K)}_s \cap A_s} \sum_{z \in \Z^{2N}_d} 
\PR [(y_i,0) \to (z_i,s-r) \; \mbox{for $i \leq N$}, \max_{i \leq K} \sigma_i \leq s-r] 
\psi^{(N)}_r(z,x) \\
& \leq & \sum_{y \in B^{(N)}_s} \sum_{y' \in B^{(K)}_s \cap A_s}\sum_{z \in \Z^{2N}_d} 
\PR  [(y_i,0) \to (z_i,s-r) \; \mbox{for $i \leq N$}, \max_{i \leq K} \sigma_i \leq s-r]  \\
&& \hspace{2in}
\left(c_9(x) (\log r)^{-{N \choose 2}}\prod_{i=1}^N p_r(z_i) + c_{10}(x,\mu) r^{-N} \log^{-\mu}\! r \right) \\
& = & I + II.
\end{eqnarray*}
The second term $II$ can be bounded by
\begin{eqnarray*}
 C r^{-N} \log^{-\mu}\! r  \sum_{y \in B^{(N)}_s} \sum_{y' \in B^{(K)}_s}\sum_{z \in \Z^{2N}_d} 
\PR_{y}[X_{s-r}=z] 
& \leq & C r^{-N} \log^{-\mu}\! r |B^{(N)}_s| \, |B^{(K)}_s| \\
& \leq & C t^K (\log t)^{\eta N - (\frac12-2\epsilon) (N+K) - \mu}
\end{eqnarray*}
and (after $\eta$ is chosen above) we may choose $\mu$ large enough that this 
does not affect the asymptotics. 

For the term $I$ we perform the sum over $z$ and bound it as in (\ref{pairbound}) by 
\begin{eqnarray*}
&& \hspace{-.4in} C(x) (\log r)^{-{N \choose 2}} 
 \sum_{y \in B^{(N)}_s} \sum_{y' \in B^{(K)}_s \cap A_s} 
 \PR_{y,y'}[X^{(i)}_s = x_i\; \mbox{for $i \leq N$}, \max_{j \leq K} \sigma_j \leq s-r] \\
 & = & C(x) (\log r)^{-{N \choose 2}} 
 \sum_{y \in B^{(N)}_s} \sum_{y' \in B^{(K)}_s \cap A_s} \prod_{i=1}^K P_{(y_i,y'_i)}[ X^{(1)}_s = x_i, \tau_{12} \leq s].
 \prod_{i=K+1}^N P_{y_i}[ X_s = x_i]
\end{eqnarray*}
For $x \in \{x_1,\ldots,x_K\}$ we use 
\begin{eqnarray*}
&& \hspace{-.5in} \sum_{y \in B^{(1)}_s} \sum_{y' \in B^{(1)}_s} 
P_{(y,y')}[ X^{(1)}_s = x, \tau_{12} \leq s] \,\I(|y'-y| > s^{1/2} \log^{-\gamma}\!s) \\
&= & \sum_{y \in B^{(1)}_s} \sum_{y^{''} \in B^{(1)}_s-y} P_{(0,y^{''})}[ X^{(1)}_s = x-y, \tau_{12} \leq s] 
\, \I(|y^{''}| > s^{1/2} \log^{-\gamma}\!s) \\
& \leq & \sum_{y \in \Z^2} \sum_{y^{''} \in B^{(1)}_{s}-B^{(1)}_{s}}
 P_{(0,y^{''})}[ X^{(1)}_s = x-y, \tau_{12} \leq s] \, \I(|y^{''}| > s^{1/2} \log^{-\gamma}\!s) \\
& = & \sum_{z \in B^{(1)}_{s}-B^{(1)}_{s}} P_{z}[\tau_0 \leq 2s] \, \I(|z| > s^{1/2} \log^{-\gamma}\!s) \\
& \leq & C  \frac{\log \log s}{\log s} |B^{(1)}_{s}-B^{(1)}_{s}|  \quad \mbox{by Lemma \ref{Gbounds} (a)} \\
& \leq & C  \frac{s\log \log s}{\log^{1-2 \epsilon}\!s}.
\end{eqnarray*}
Then we bound the term $I$ by 
\[
 C (\log r)^{-{N \choose 2}} \left(\frac{s\log \log s}{\log^{1-2 \epsilon}\!s}\right)^K
= O\left( (\log t)^{-{N \choose 2}} (t \log \log t)^K (\log t)^{-(\frac12 - 2 \epsilon) K} \right)
\]
which completes step 2 and hence the correction term to the asymptotics in (\ref{MR2}) (where we 
choose $\delta = 2 \epsilon N$).

\noindent
\textbf{Proof of claims (\ref{claim2000a}) and (\ref{claim2000b}).} 
For $x \in \{x_1,\ldots,x_N\}$
\begin{eqnarray*}
\sum_{y \in \Z^2} \sum_{y' \in \Z^2} \PR_{(y,y')}[ X^{(1)}_s=x, \tau_{12} \leq s] 
& = & \sum_{y \in \Z^2} \sum_{y^{''} \in \Z^2} \PR_{(0,y^{''})}[ X^{(1)}_s=x+y, \tau_{12} \leq s]  \\
& = & \sum_{y^{''} \in \Z^2} \PR_{(0,y^{''})}[\tau_{12} \leq s] \\
& = & \sum_{z \in \Z^2} \PR_{z}[\tau_{0} \leq 2s] \\
& = & \sum_{z \in \Z^2} \PR_{0}[\tau_{z} \leq 2s].
\end{eqnarray*}
The range $\mathcal{R}_{2s} = \sum_z \I(\tau_z \leq 2s)$ of a simple random walk over the interval $[0,2s]$ is 
bounded by the number of steps of the walk - a Poisson ($2s$) variable - and this leads to (\ref{claim2000a}).
Similarly, 
\begin{eqnarray*}
&& \hspace{-.4in} \sum_{y \in \Z^2} \sum_{y' \in \Z^2} \PR_{(y,y')}[ X^{(1)}_s=x, \tau_{12} \leq s] 
\I(|y'| > s^{1/2} \log^{\epsilon}\!s)
\\ & = & \sum_{y \in \Z^2} \sum_{y^{''} \in \Z^2} \PR_{(0,y^{''})}[ X^{(1)}_s=x+y, \tau_{12} \leq s] 
\I(|y''-y| > s^{1/2} \log^{\epsilon}\!s) \\
& \leq & \sum_{y \in \Z^2} \sum_{y^{''} \in \Z^2} \PR_{(0,y^{''})}[ X^{(1)}_s=x+y, \tau_{12} \leq s] \\
&& \hspace{1in}
\left( \I(|y^{''}|  > \frac{s^{1/2} \log^{\epsilon}\!s}{2}) + \I(|y| > \frac{s^{1/2} \log^{\epsilon}\!s}{2} > |y^{''}|)  \right) \\
& \leq & \sum_{y^{''} \in \Z^2} \PR_{(0,y^{''})}[\tau_{12} \leq s] \I(|y^{''}|  > \frac{s^{1/2} \log^{\epsilon}\!s}{2}) \\
&& \hspace{.3in} + \sum_{y \in \Z^2} \sum_{y^{''} \in \Z^2} \PR_{(0)}[ X^{(1)}_s=x+y] \I(|y| > \frac{s^{1/2} \log^{\epsilon}\!s}{2} > |y^{''}|) \\
& \leq & \sum_{z \in \Z^2} \PR_{0} [\tau_{z} \leq 2s] \I(|z|  > \frac{s^{1/2} \log^{\epsilon}\!s}{2}) + C s \log^{2 \epsilon}\!s \;
 \PR_{(0)}\left[ |X^{(1)}_s-x| \geq  \frac{s^{1/2} \log^{\epsilon}\!s}{2}\right] \\
 & = & I +II.
\end{eqnarray*}
The second term $II$ is $O(s \log^{-\mu}\!s)$ by a large deviation estimate. The first term $I$ can be bounded 
by 
\[
I = \E_0\left[ |\mathcal{R}_{2s} \cap (\frac12 B^{(1)}_s)^c|\right]  \leq \E_0[ \mathcal{R}_{2s}^2]^{1/2}
\PR_0[ \mathcal{R}_{2s} \cap (\frac12 B^{(1)}_s)^c \neq \emptyset]^{1/2}
\]
by the Cauchy Schwarz inequality. The expectation is $\E_0[ \mathcal{R}_{2s}^2] = 2s +4s^2$,
while the probability $\PR_0[ \mathcal{R}_{2s} \cap (\frac12 B^{(1)}_s)^c \neq \emptyset]$
is  $O(\log^{-2 \mu}\!s)$ again by the large deviation estimate (\ref{LD}). This shows that the term $I$ is
also $O(s \log^{-\mu}\!s)$ completing the proof of  (\ref{claim2000b}).

\noindent
\textbf{Proof of Lemma \ref{Jamieplus}.}
For the bound (\ref{JLplus}) we may follow closely the argument of Lemma 12 in \cite{KvdB1}.
It follows again a series of approximations based on splitting at time $s = t \log^{-\alpha}\!t$ for $\alpha > 2$.
\begin{eqnarray*}
&& \hspace{-.3in} \sum_{y} \left|\psi^{(N)}_t(x,y) - c_0(x) (\log t)^{-{N \choose 2}} p_t(y) \right| \\
& = & \sum_{y} \left|\PR_x[X_t = y, \tau_c >t] - c_0(x) (\log t)^{-{N \choose 2}} p_t(y) \right|  \\
& = & \sum_{y} \left|\PR_x[X_t = y, \tau_c >t] - c_0(x) (\log s)^{-{N \choose 2}} p_t(y)\right| +  \mathcal{E}_1 \\
& = &  \sum_{y} \left|\PR_x[X_t = y, \tau_c >s] - c_0(x) (\log s)^{-{N \choose 2}} p_t(y) \right| + \mathcal{E}_2 \\
& = &  \sum_{y} \left| \sum_z \left(\PR_x[X_s=z, \tau_c >s] - c_0(x) (\log s)^{-{N \choose 2}} p_s(z)\right)
p_{t-s}(y-z) \right| + \mathcal{E}_2\\
& = & \sum_{y} \left| \sum_z \left(\PR_x[X_s=z, \tau_c >s] - c_0(x) (\log s)^{-{N \choose 2}} p_s(z)\right)
p_{t-s}(y) \right| + \mathcal{E}_3\\
& = & \sum_{y} \left|  \left(\PR_x[\tau_c >s] - c_0(x) (\log s)^{-{N \choose 2}} \right)
p_{t-s}(y) \right| + \mathcal{E}_3 \\
& \leq & \left|\PR_x[\tau_c >s] - c_0(x) (\log s)^{-{N \choose 2}} \right| + \mathcal{E}_3.
\end{eqnarray*}
The error bound (\ref{MR1}) for the asymptotics for $\p(t)$ complete the proof if we can show that
$\mathcal{E}_3(t) = O \left((\log t)^{-{N \choose 2}-1+\epsilon} \right)$. 

The error $\mathcal{E}_1 = O((\log t)^{-{N \choose 2}-1}\log \log t)$ due to the choice of $s$. The second error
$D \mathcal{E}_2 = \mathcal{E}_2 - \mathcal{E}_{1}$ is precisely $P_x[\tau_c \in (s,t]] = \p(t)-\p(s)$
and the asymptotics (\ref{MR1}) show that this is of the same order as $\mathcal{E}_1$.
To control $D \mathcal{E}_3$ we may replace first the sum all $z$ by $z \in B_s^{(N)}$ (as in (\ref{Bbox})), and then the sum over all $y$ by $y \in B_{t}^{(N)}$, all at the expense of an error that is $O(\log^{-\mu}\!t)$ for any $\mu$ (we omit the details which are similar to previous steps). This leaves
\[
D \mathcal{E}_3 \leq \sum_{y \in B_{t}^{(N)}} \sum_{z \in B_s^{(N)}} \left(\PR_x[X_s=z, \tau_c >s] + c_0(x) (\log s)^{-{N \choose 2}} p_s(z)\right)
|p_{t-s}(y-z) - p_{t-s}(y)|. 
\]
Repeated uses of the one dimensional local central limit theorem (\ref{lclt}) lead to the bound
\[
|p_{t-s}(y-z) - p_{t-s}(y)| \leq C|z| t^{-N-\frac12} + O(t^{-N+1-\nu}) \quad \mbox{for any $\nu < 3/2$}
\]
and then 
\begin{eqnarray*}
D \mathcal{E}_3 & \leq & C (\log s)^{-{N \choose 2}} \sum_{y \in B_{t}^{(N)}}\left(t^{-N-\frac12} \max_{z \in B_s^{(N)}} |z|  + O(t^{-N+1-\nu}) \right) \\
& \leq & C (\log s)^{-{N \choose 2}} (t \log^{2 \epsilon}\!t)^N \left(t^{-N-\frac12} s^{1/2} \log^{\epsilon}\!s  + O(t^{-N+1-\nu}) \right) \\
& = & O \left((\log t)^{-{N \choose 2}-1}\right) 
\end{eqnarray*}
by the choice of $s = t \log^{-\alpha}\!t$ with $\alpha > 2$ large enough. 

To prove (\ref{UBNplus}) we argue as in Lemma \ref{Gbounds} (b): we 
split $[0,t]$ at a time $s = t \log^{-\alpha}\!t$ for some $\alpha>0$. 
\begin{eqnarray*}
\psi^{(N)}_t(x,y) &=&\sum_z  \psi^{(N)}_s(x,z)  \psi^{(N)}_{t-s}(z,y) \\
& \leq & \sum_z \psi^{(N)}_s(x,z)  p_{t-s}(y-z) \\
& = &  \sum_z \psi^{(N)}_s(x,z) p_t(y) + \mathcal{E}_5 \\
& = &  \PR_{x}[\tau_c >s] p_t(y) + \mathcal{E}_5.
\end{eqnarray*}
It remains to show that $\mathcal{E}_5 = O(t^{-N} \log^{-\mu}\!t)$.
We may again restrict to summing over $|z| \in B^{(N)}_s$ and then 
estimates on the derivatives of the Gaussian transition density show
\begin{eqnarray*}
|p_{t-s}(z-y) - p_t(y)| & \leq & |p_{t-s}(z-y) - p_t(z-y)| + |p_{t}(z-y) - p_t(y)| \\
& \leq & C s t^{-N-1} + C |z| t^{-N-\frac12} + O(t^{-N+1-\nu})
\end{eqnarray*}
which again leads to the desired bound by taking $\alpha = \alpha(\mu)$ large enough. 

\noindent
Acknowledgments. Roger Tribe is partially supported by a Leverhulme Research Fellowship RF-2-16-655.
Jamie Lukins is supported by EPSRC as part of the MASDOC DTC at the University of Warwick. Grant No. EP/HO23364/1.

%
%
\end{document}